\input ./style/arxiv-general.cfg
\documentclass[MSNbibl,number,citesort,seceqn,dvips]{arxbj}
\makeatletter
   \@ifpackageloaded{graphicx}{}{\usepackage{graphicx}}
\makeatother
\volume{22}
\issue{4}
\pubyear{2016}
\firstpage{2372}
\lastpage{2400}
\doi{10.3150/15-BEJ732}% Updated by VTEXPTS2LaTeX.exe, 17.08.2015 08:13
\docsubty{FLA}

\makeatletter
\def\eqref#1{(\ref{#1})}

\newcommand{\rrvert}{\vert}

\newcommand{\rrVert}{\Vert}
\newcommand{\llvert}{\vert}
\newcommand{\llVert}{\Vert}
\newtheorem{theorem}{Theorem}[section]
\newtheorem{lemma}{Lemma}
\newremark{remarks}{Remarks}
\makeatother

\begin{document}
\begin{frontmatter}

\title{Asymptotic theory for statistics of the Poisson--Voronoi approximation}
\runtitle{Statistics of the Poisson--Voronoi approximation}

\begin{aug}
%%%% inicialai - be tarpu
% Corresponding author: Joseph Yukich - jey0@lehigh.edu% Updated by
%VTEXPTS2LaTeX.exe, 17.08.2015 08:13
\author[A]{\inits{C.}\fnms{Christoph}~\snm{Th\"ale}\thanksref{A}\ead[label=e1]{christoph.thaele@rub.de}}%,
%\author[]{\inits{}\fnms{}~\snm{}\thanksref{}\ead[label=]{}}
\and
\author[B]{\inits{J.E.}\fnms{J.E.}~\snm{Yukich}\corref{}\thanksref{B}\ead[label=e2]{joseph.yukich@lehigh.edu}}
%%\runauthor{} %% auto
%\dedicated{}
\address[A]{Faculty of Mathematics, Ruhr University Bochum, Bochum,
Germany.\\ \printead{e1}}
\address[B]{Department of Mathematics, Lehigh University, Bethlehem,
PA 18015, USA.\\ \printead{e2}}
\end{aug}

% HISTORY:
%
\received{\smonth{12} \syear{2014}}% Updated by VTEXPTS2LaTeX.exe,
%17.08.2015 08:13
%
\revised{\smonth{4} \syear{2015}}% Updated by VTEXPTS2LaTeX.exe,
%17.08.2015 08:13

% ABSTRACT
%
\begin{abstract}
This paper establishes expectation and variance asymptotics for
statistics of the Poisson--Voronoi approximation of
general sets, as the underlying intensity of the Poisson point process
tends to infinity. Statistics of interest include volume,
surface area, Hausdorff measure, and the number of faces of
lower-dimensional skeletons. We also consider the complexity of the
so-called Voronoi zone and the iterated Voronoi approximation. Our
results are consequences of general limit theorems proved with an
abstract Steiner-type formula applicable in the setting of sums of
stabilizing functionals.
\end{abstract}

% KEYWORDS
% visi is mazosios raides ir pagal abecele
%
\end{frontmatter}
%
%s1
\section{Main results}\label{sect1}

The Poisson--Voronoi mosaic is a classical and prominent example of a
random mosaic and is used in a wide range of fields, including
astronomy, biology, material sciences and telecommunications.
%with applications ranging from biology and material sciences to recent
%developments in telecommunication.
If ${\mathcal{P}}_\lambda $ is a Poisson point process on
$Q:=[-1/2,1/2]^d$ whose
intensity measure has density $\lambda \kappa ( \cdot )$ with
respect to
the Lebesgue measure ($d\geq2$, $\lambda \in(0, \infty)$ and
$\kappa $ is a
continuous function on $Q$ bounded away from zero and infinity), the
Voronoi cell $v(x):=v(x, {\mathcal{P}}_\lambda )$ associated with
$x\in{\mathcal{P}}_\lambda $ is
the set of all $z\in Q$ such that the distance between $z$ and $x$ is
less than the distance between $z$ and any other point of ${\mathcal
{P}}_\lambda $.
Clearly, $v(x)$ is a random convex polytope and the collection of all
$v(x)$ with $x\in{\mathcal{P}}_\lambda $ partitions $Q$ and is
called the
Poisson--Voronoi mosaic of $Q$.

Let $A\subset Q$ be a full-dimensional admissible set %in the sense of
%Section~\ref{sec:GeneralTheorems} below
whose boundary has positive and finite \mbox{$(d-1)$}-dimensional Hausdorff measure.
Admissible sets, formally defined in Section~\ref{sec:GeneralTheorems}, include in particular, convex sets, sets of
positive reach, differentiable manifolds with smooth boundary as well
as certain finite unions of such sets. Given such $A \subset Q$, the
Poisson--Voronoi approximation $\operatorname{PV}_\lambda (A)$ of $A$ is the union
of all Voronoi cells $v(x)$ with $x\in A$, that is,
\[
\operatorname{PV}_\lambda (A):=\bigcup_{x\in{\mathcal{P}}_\lambda \cap
A}v(x).
\]
Typically, $A$ is an unknown set having unknown geometric
characteristics such as volume or surface area. Notice that $\operatorname{PV}_\lambda (A)$ is a random polyhedral approximation of $A$, which closely
approximates $A$ as $\lambda $ becomes large. One might expect that the
volume and surface area of $\operatorname{PV}_\lambda (A)$, respectively
denoted by
$V_\lambda (A)$ and $S_\lambda (A)$, also closely approximate the
volume and
surface area of $A$. Our first goal is to show that this is indeed the
case, though the surface area asymptotics involve a universal
correction factor, denoted by $c_2$ in the sequel.
%which closely approximates $A$ for, which is getting finer and finer,
%as $\la$ tends to infinity. By $V_\la(A)$ and $S_\la(A)$ we denote the
%volume and the surface area of $\operatorname{PV}_\la(A)$.
%Our first goal is to prove expectation and variance % asymptotics for
%$V_\la(A)$ and $S_\la(A)$.
For sets $A$ which are convex or which have a smooth boundary,
first-order asymptotics have been previously established in \cite
{HevelingReitzner,Penrose07b,RSZ,Yu};
second-order asymptotics for sets $A$ having a smooth boundary are
given in \cite{Yu}, while \cite{SchulteVoronoi} provides second-order
inequalities when $A$ is a convex set. %One of the main
%achievements of this paper is that we can deal with a general class of
%sets $A$, %comprising these special cases.
We extend the limit theory of these papers and obtain first- and
second-order asymptotics whenever $A$ belongs to the more general class
of admissible sets. In particular, we show that the variance
asymptotics for $V_\lambda (A)$ are proportional to the $\kappa $-weighted
surface content of $A$, resolving a conjecture implicit in Remark~2.2
of \cite{RSZ}. The approach relies on a general and far-reaching
Steiner-type formula from \cite{HLW}, together with stabilization
properties of geometric functionals of the Poisson--Voronoi
mosaic.

In the sequel, we write $f(\lambda )\sim c g(\lambda )$ for real-valued
functions $f$ and $g$ and constants \mbox{$c \in[0, \infty)$} if $\lim_{\lambda \to\infty} f(\lambda )/g(\lambda ) = c$.
Throughout, we denote
the $s$-dimensional Hausdorff measure by ${\mathcal{H}}^s$, $s \in[0,
\infty
)$. Furthermore, we say that $\partial A$ contains a subset $\Gamma$
of differentiability class $C^2$ with ${\mathcal{H}}^{d-1}(\Gamma)\in
(0,\infty
)$ if $\Gamma\subset\partial A$ is an open and twice differentiable
$(d-1)$-dimensional sub-manifold in $\mathbb{R}^d$ in the usual sense of
differential geometry. Finally, for $\gamma\in\mathbb{R}$ we define the
{\em$\kappa $-weighted surface content}\vspace*{-2pt}
\[
{\mathcal{H}}_{\kappa,\gamma}^{d-1}(\partial A):=\int
_{\partial
A}\kappa (x)^{1-\gamma/d} {\mathcal{H}}^{d-1}(
\mathrm{d}x).
\]
Observe that ${\mathcal{H}}_{\kappa,\gamma}^{d-1}(\partial A)$
reduces to the
usual surface content ${\mathcal{H}}^{d-1}(\partial A)$ of $\partial
A$ if either
$\gamma= d$ and $\kappa $ is arbitrary or $\kappa\equiv1$ and
$\gamma
\in\mathbb{R}$ is arbitrary.
%for arbitrary $\ka$, $\gamma=d$ or if $\kappa\equiv1$ for arbitrary $
%\gamma\in\RR$.

\begin{theorem}\label{thm:VoronoiV+S}
There are constants $c_1,c_2 \in(0, \infty)$ depending only on the
dimension $d$ such that\vspace*{-1pt}
\[
\mathbb{E}V_\lambda (A) - V(A)\sim c_1 \lambda
^{-{1/d}} {\mathcal{H}}^{d-1}(\partial A) \quad \mbox{and}\quad
\mathbb{E}S_\lambda (A) \sim c_2 {\mathcal
{H}}_{\kappa,d-1}^{d-1}(\partial A).
\]
Moreover, there are constants $c_3,c_{4,1},c_{4,2} \in[0, \infty)$
depending only on $d$ such that\vspace*{-1pt}
\[
\operatorname{Var} \bigl[V_\lambda (A) \bigr]\sim c_3
\lambda ^{-1-{1/d}} {\mathcal{H}}^{d-1}(\partial A)
\]
and\vspace*{-1pt}
\[
\operatorname{Var} \bigl[S_\lambda (A) \bigr]\sim\lambda ^{-1+{1/d}}
\bigl(c_{4,1}{\mathcal{H}}_{\kappa,2(d-1)}^{d-1}(\partial
A)+c_{4,2}{\mathcal{H}}_{\kappa
^2,d-1}^{d-1}(\partial A)
\bigr).
\]
If $\partial A$ contains a subset $\Gamma$ of differentiability class
$C^2$ with $\mathcal{H}^{d-1}(\Gamma) \in(0, \infty)$, and if
$\kappa
\equiv1$, then $c_3$ and $c_4:=c_{4,1}+c_{4,2}$ are strictly positive.
\end{theorem}

Next, we turn to other metric parameters of the Poisson--Voronoi
approximation, which can be handled by our general set-up. To this end,
for $\ell\in\{0,\ldots,d-1\}$ denote by $\operatorname{skel}_\ell(\operatorname{PV}_\lambda (A))$ the union of all $\ell$-dimensional faces belonging to
$\partial(\operatorname{PV}_\lambda (A))$, the boundary of $\operatorname{PV}_\lambda
(A)$, and
let $H^{(\ell)}_\lambda (A)$ be the $\ell$-dimensional Hausdorff measure
of $\operatorname{skel}_\ell(\operatorname{PV}_\lambda (A))$. More formally, if
$\mathcal{F}_\ell(P)$ stands for the collection of $\ell
$-dimensional faces of
a polytope $P$, then\vspace*{-1pt}
\[
H^{(\ell)}_\lambda(A):=  {1 \over  d - \ell} \mathop{\sum_{x\in\mathcal{P}_\lambda}}_{x\in A}
\mathop{\sum_{f\in{\mathcal{F}_\ell(v(x))}}}_{f\subset\partial({\operatorname{PV}}_\lambda(A))}{\mathcal{H}}^{\ell}(f).
\]
%\[
%H^{(\ell)}_\lambda (A):=\mathop{\sum_{x\in{\mathcal{P}}_\lambda}}_{x\in A}
%\mathop{\sum_{f\in
%\mathcal{F}_\ell(v(x))}}_{ f\subset\partial(\operatorname{PV}_\lambda
%(A))}{\mathcal{H}} ^{(\ell)}(f).
%\]
%
Note that $H_\lambda^{(d-1)}(A)$ coincides with $S_\lambda(A)$
considered in Theorem~\ref{thm:VoronoiV+S}.

\begin{theorem}\label{thm:VoronoiSkeleton}
Let $\ell\in\{0,\ldots,d-1\}$. Then there are constants $c_5 \in(0,
\infty)$ and $c_{6,1},c_{6,2} \in[0, \infty)$ depending only on $d$
and $\ell$ such that
\[
\mathbb{E}H^{(\ell)}_\lambda (A) \sim c_5 \lambda
^{1 - {1/d} - {\ell/d} } {\mathcal {H}}_{\kappa,\ell
}^{d-1}({\partial A})
\]
and
\[
\operatorname{Var} \bigl[H^{(\ell)}_\lambda (A) \bigr] \sim\lambda
^{1-{1/
d} -{(2 \ell)/d}} \bigl(c_{6,1} {\mathcal{H}}_{\kappa,2\ell}^{d-1}({
\partial A})+c_{6,2} {\mathcal {H}}^{d-1}_{\kappa^2,\ell
}(
\partial A) \bigr).
\]
If $\partial A$ contains a subset $\Gamma$ of differentiability class
$C^2$ with $\mathcal{H}^{d-1}(\Gamma) \in(0, \infty)$, and if
$\kappa
\equiv1$, then $c_6:=c_{6,1}+c_{6,2}$ is strictly positive.
\end{theorem}

With the exception of $H_\lambda^{(0)}(A)$, the number of vertices on
$\partial(\operatorname{PV}_\lambda(A))$, we have investigated only metric
parameters of the Poisson--Voronoi approximation, namely the volume,
the surface area and the Hausdorff measure of lower-dimensional
skeletons. On the other hand, the combinatorial complexity of $\operatorname{PV}_\lambda (A)$ is also of interest. For example, it is natural to
ask how
many %vertices or facets, or more generally, how many
$\ell$-dimensional faces ($\ell\in\{0,\ldots,d-1\}$) belong to
$\partial(\operatorname{PV}_\lambda(A))$. In contrast to volume and surface
area, combinatorial parameters of the Poisson--Voronoi approximation
have apparently not been studied in the literature. The general theory
developed in Section~\ref{sec:GeneralTheorems} allows us to
investigate such parameters. To state the result, for $\ell\in\{
0,\ldots,d-1\}$ we let $f_\lambda ^{(\ell)}(A)$ be the number of
$\ell
$-dimensional faces belonging to $\partial(\operatorname{PV}_\lambda (A))$. Note
that $f_\lambda^{(0)}(A)=H_\lambda^{(0)}(A)$. % investigated in
%Theorem~\ref{thm:VoronoiSkeleton}.

\begin{theorem}\label{thm:VoronoiCombinatorics}
Let $\ell\in\{0,\ldots,d-1\}$. Then there are constants $c_7 \in(0,
\infty)$ and $c_{8,1},c_{8,2} \in[0, \infty)$ depending only on the
dimension $d$ and on $\ell$ such that
\[
\mathbb{E}f_\lambda ^{(\ell)}(A) \sim c_7 \lambda
^{1-{1/d}} {\mathcal{H}}_{\kappa,0}^{d-1}(\partial A)
\]
and
\[
\operatorname{Var} \bigl[f_\lambda ^{(\ell)}(A) \bigr]\sim\lambda
^{1-{1/
d}} \bigl(c_{8,1} {\mathcal{H}} _{\kappa,0}^{d-1}(
\partial A)+c_{8,2} {\mathcal{H}}_{\kappa
^2,0}^{d-1}(
\partial A) \bigr).
\]
If $\partial A$ contains a subset $\Gamma$ of differentiability class
$C^2$ with $\mathcal{H}^{d-1}(\Gamma) \in(0, \infty)$, and if
$\kappa
\equiv1$, then $c_8:=c_{8,1}+c_{8,2}$ is strictly positive.
\end{theorem}

%A new aspect arises when we consider the
Next, we consider certain functionals of Voronoi cells intersecting
only a part of the boundary of $A$. Formally,
given an admissible set $A$ and $A_0\subset\partial A$ such that
${\mathcal{H}}
^{d-1}(A_0) \in(0, \infty)$, define the Poissson--Voronoi zone $\operatorname{PVZ}_\lambda (A_0)$ of $A_0$ by
\[
\operatorname{PVZ}_\lambda (A_0):=\mathop{\bigcup
_{x\in{\mathcal{P}}_\lambda }}_{
v(x)\cap A_0\neq
\varnothing}v(x).
\]
Given $\ell\in\{0,\ldots,d-1\}$, let $\widehat{f}_\lambda ^{(\ell
)}(A_0)$ denote the number of $\ell$-dimensional faces of $\operatorname{PVZ}_\lambda (A_0)$. We emphasize that this construction is very
similar to
the construction of a zone in a hyperplane arrangement; see \cite
{Matousek}. Following these classical ideas, we define the complexity
of $\operatorname{PVZ}_\lambda (A_0)$ as $\operatorname{Co}_\lambda (A_0):=\widehat
{f}_\lambda
^{(0)}(A_0)+\cdots+\widehat{f}_\lambda ^{(d-1)}(A_0)$. The zone theorem
in discrete geometry (see Theorem~6.4.1 in \cite{Matousek}) asserts
that the complexity of a zone of an arbitrary hyperplane arrangement is
of surface-order. Our next result shows a similar surface-order
behaviour for the expectation and the variance in case of a random
Poisson--Voronoi zone.

\begin{theorem}\label{thm:VoronoiZone}
There are constants $c_9 \in(0, \infty)$ and $c_{10,1},c_{10,2} \in
[0, \infty)$ depending only on $d$ such that
\[
\mathbb{E} \operatorname{Co}_\lambda (A_0)\sim c_9
\lambda ^{1-{1/d}} {\mathcal{H}}_{\kappa,0}^{d-1}(A_0)
\]
and
\[
\operatorname{Var} \bigl[\operatorname{Co}_\lambda (A_0) \bigr]\sim
\lambda ^{1-{1/d}} \bigl(c_{10,1} {\mathcal{H}} _{\kappa,0}^{d-1}(A_0)+c_{10,2}
{\mathcal{H}}_{\kappa
^2,0}^{d-1}(A_0) \bigr).
\]
If $A_0$ contains a subset $\Gamma$ of differentiability class $C^2$
with $\mathcal{H}^{d-1}(\Gamma) \in(0, \infty)$, and if $\kappa
\equiv1$,
then $c_{10}:=c_{10,1}+c_{10,2}$ is strictly positive.
\end{theorem}

Another application of our results concerns the iterated
Poisson--Voronoi approximation, defined recursively as follows:
\[
\operatorname{PV}_\lambda ^{(1)}(A):=\operatorname{PV}_\lambda (A) \quad \mbox
{and}\quad  \operatorname{PV}_\lambda ^{(n)}(A):=\operatorname{PV}_{n\lambda }
\bigl(\operatorname{PV}_\lambda ^{(n-1)}(A) \bigr)
\]
for integers $n\geq1$ (note that the intensity used in the $n$th
iteration is $n\lambda$, where $\lambda >0$ is fixed). By $V_\lambda ^{(n)}$,
$S_\lambda ^{(n)}$ and $f_\lambda ^{\ell,(n)}$ we denote the volume, the
surface area and the number of
$\ell$-dimensional faces ($\ell\in\{0,\ldots,d-1\}$) of the $n$th
iterated Poisson--Voronoi approximation, respectively. Moreover, by
$H_{\lambda}^{\ell,(n)}$ we indicate the $\ell$-dimensional
Hausdorff measure of the $\ell$-skeleton of $\operatorname{PV}_\lambda ^{(n)}(A)$,
$\ell\in\{0,\ldots,d\}$. Note that our construction of the iterated
Poisson--Voronoi approximation is close to that of so-called aggregate
mosaics introduced in %by Tchoumatchenko and Zuyev in
\cite{TchZuyev}. The expectation analysis of functionals of the
iterated Poisson--Voronoi mosaic yields the following result. Variance
asymptotics are less tractable and we shall omit them. For simplicity,
we shall assume that the Poisson point process $\mathcal{P}_\lambda $ is
homogeneous with $\kappa \equiv1$.
%\CT[In the last sentence, I changed variance analysis to variance
%asymptotics.]

%\CT[A $1$ was missing in the constant $c_{2,n}$.]
%
\begin{theorem}\label{thm:IteratedVoronoi}
Suppose that $\kappa \equiv1$ and let $c_1$ and $c_2$ be the constants
from Theorem~\ref{thm:VoronoiV+S}, $c_5$ the constant from Theorem~\ref{thm:VoronoiSkeleton}, and $c_7$ the constant from Theorem~\ref
{thm:VoronoiCombinatorics}. Put $c_{2,n}:=1+c_2+c_2^2+\cdots
+c_2^{n-1}$ for integers $n\geq1$. Then%\vspace*{9pt}
\begin{eqnarray*}
\mathbb{E}V_\lambda ^{(n)}-V(A) &\sim& c_{1}
c_{2,n} \lambda ^{-{1/d}} {\mathcal{H}} ^{d-1}({\partial
A}),\vphantom{\sum^.}
\\
\mathbb{E}S_\lambda ^{(n)}- S(A) & \sim& c_{2}
c_{2,n} \lambda ^{-{1/d}} {\mathcal{H}} _{\kappa,d-1}^{d-1}({
\partial A}),\vphantom{\sum^.}
\\
\mathbb{E}H_\lambda ^{\ell,(n)} &\sim& c_5
c_{2,n} \lambda ^{1-{1/d}-{\ell
/ d}} {\mathcal{H}}_{\kappa,\ell}^{d-1}({
\partial A}),\vphantom{\sum^.}
\\
\mathbb{E}f_\lambda ^{\ell,(n)} &\sim& c_{7}
c_{2,n} \lambda ^{1-{1/d}} {\mathcal{H}} _{\kappa,0}^{d-1}({
\partial A}).\vphantom{\sum^.}
\end{eqnarray*}
\end{theorem}

\begin{remarks*}
\begin{longlist}%[(iii)]
\item[(i)] \textit{Theorem~\ref{thm:VoronoiV+S} (related work).} The
set $\operatorname{PV}_\lambda (A)$ was introduced in \cite{KT} where it was shown
that $\lim_{\lambda \to\infty} \operatorname{Vol}(A \Delta
A_\lambda ) = 0$
almost surely, but only when $d = 1$. This almost sure limit was
extended in \cite{Penrose07b} to all dimensions $d \geq1$. When
${\mathcal{P}}
_\lambda $ denotes a homogeneous Poisson point process on $\mathbb
{R}^d$ having
intensity $\lambda $, we have that $V_\lambda (A)$ is an unbiased
estimator of
$V(A)$ (cf. \cite{RSZ}), which makes $\operatorname{PV}_\lambda (A)$ of interest
in image analysis, non-parametric statistics and quantization;
see also Section~1 of \cite{KT} and Section~1 of \cite{HevelingReitzner}.
\item[(ii)] \textit{Invariance of limits with respect to geometry.}
The common thread linking our results is that the first- and
second-order asymptotic behaviour of our functionals are geometry
independent. By this we mean that the mean and variance asymptotics are
not influenced by the precise geometric structure of the given
admissible set $A$, but
are rather controlled only by the $\kappa $-weighted surface content
of $A$.
\item[(iii)] \textit{The constants in Theorems \ref
{thm:VoronoiV+S}--\ref{thm:IteratedVoronoi}.} The explicit dependency
of the constants $c_i, i \geq1$, in Theorems \ref
{thm:VoronoiV+S}--\ref{thm:IteratedVoronoi} on the dimension $d$ and
the parameter $\ell$ is given explicitly in the general results of
Section~\ref{sec:GeneralTheorems}, especially the upcoming limits \eqref{eq: newexpect} and
\eqref{eq: newvar}. More precisely, let $\mathcal{P}_1^{\mathrm{hom}}$
be a
homogeneous Poisson point process on $\mathbb{R}^d$ of unit intensity
and put $\mathbb{R}_+^{d-1}:= \mathbb{R}^{d-1} \times\mathbb{R}^+$. Let
\[
\operatorname{PV} \bigl(\mathbb{R}^{d-1}_+ \bigr):= \bigcup
_{x \in\mathcal{P}_1^{\mathrm{hom}} \cap\mathbb{R}
_+^{d-1}} v(x)
\]
be the Poisson--Voronoi approximation of $\mathbb{R}^{d-1}_+$. Then the
general results show that the
expectation and variance asymptotics are controlled by the $\kappa
$-weighted surface content of $A$ as well as by
the expected behaviour of metric and combinatorial parameters of the
simpler object $\operatorname{PV}(\mathbb{R}^{d-1}_+)$.
Finding explicit numerical values for the constants $c_i, i \geq1$,
arising in expectation and variance asymptotics is a separate problem
which we do not tackle here.
\item[(iv)] \textit{Extensions of Theorems \ref
{thm:VoronoiV+S}--\ref{thm:IteratedVoronoi}.} By Theorem~\ref
{eq:WLLN} below, the expectation asymptotics in Theorems \ref
{thm:VoronoiV+S}--\ref{thm:IteratedVoronoi} may be upgraded to a weak
law of large numbers holding in the $L^1$- and $L^2$-sense.
\item[(v)] \textit{General surface-order results.} Although Theorems
\ref{thm:VoronoiV+S}--\ref{thm:IteratedVoronoi} only deal with
statistics of the Poisson--Voronoi approximation, we emphasize that
they follow from general theorems (presented in Section~\ref{sec:GeneralTheorems} below) for general surface-order stabilizing
functionals. These general theorems are applicable in a wider context,
establishing, for example, expectation and variance asymptotics for the
number of maximal points in a random sample, as described in Remark
(iii) after Theorem~\ref{eq:Var}.
\end{longlist}
\end{remarks*}

The rest of this paper is structured as follows. In Section~\ref{sec:GeneralTheorems}, we make precise our framework, in particular, we
introduce the class of admissible sets and score functions. We also
state there two general theorems which yield Theorems \ref
{thm:VoronoiV+S}--\ref{thm:IteratedVoronoi}. Their proofs form the
content of Section~\ref{sec:ProofsGeneralTheorems}, while Section~\ref{sec:ProofsApplications} contains the proofs of Theorems \ref
{thm:VoronoiV+S}--\ref{thm:IteratedVoronoi}. Section~\ref{sec:VLB}
establishes the asserted variance lower bounds in Theorems \ref
{thm:VoronoiV+S}--\ref{thm:VoronoiZone}.

%s2
\section{Framework and general theorems}\label{sec:GeneralTheorems}

Let ${\mathcal{P}}_\lambda $ denote a Poisson point process on
$\mathbb{R}^d$ for some $d
\geq2$ whose intensity measure has density $\lambda \kappa $ with
respect to
the Lebesgue measure on $\mathbb{R}^d$, where $\lambda \in(0, \infty
)$ but now
$\kappa $ is a bounded function on $\mathbb{R}^d$ not necessarily
bounded away
from zero. %, and $\la\in(0, \infty)$ is a %finite constant.
Furthermore, let $A \subset\mathbb{R}^d$ be a closed set such that its
boundary $\partial A$ has finite $(d-1)$-dimensional Hausdorff measure.
We consider in this section general statistics of the form
\begin{equation}
\label{eq:orig-sum}\sum_{x \in{\mathcal
{P}}_\lambda } \xi(x, {\mathcal{P}}
_\lambda, \partial A),
\end{equation}
where $\xi$ is a certain score function, which associates to a point
$x\in{\mathcal{P}}_\lambda $ a real number, which is allowed to
depend on the
surrounding point configuration ${\mathcal{P}}_\lambda $ as well as
on the set $A$
via its boundary $\partial A$. To introduce a re-scaled version and to
simplify notation, we use the abbreviation $\xi_\lambda (x, {\mathcal
{P}}_\lambda,
\partial A):=\xi(\lambda ^{1/d}x, \lambda ^{1/d}{\mathcal
{P}}_\lambda, \lambda ^{1/d}(\partial
A))$ and define
\begin{equation}
\label{eq:gsum} H^\xi({\mathcal{P}}_\lambda, \partial A):=
\sum_{x \in{\mathcal
{P}}_\lambda }\xi_\lambda (x, {\mathcal{P}}
_\lambda, \partial A).
\end{equation}
The focus of this paper is on score functions which depend on the
geometry of the set $A$ in that $\xi(x,{\mathcal{P}}_\lambda,\partial A)$ decays
with the distance of $x$ to $\partial A$. Moreover, we require $\xi$
to satisfy a weak spatial dependency condition.

To make the framework precise, we first introduce terminology, including
the collection $\mathbf{A}(d)$ of admissible sets $A\subset\mathbb
{R}^d$ as well
as the collection $\Xi$ of admissible score functions.
%\CT[Next paragraph slightly re-phrased.]
The reader may wonder about our choice of admissible sets. The
admissible sets described below have the attractive feature that their
so-called extended support measures are `well-behaved' and satisfy a
Steiner-type formula \eqref{eq:Steiner}, which is a far reaching
consequence of the classic Steiner formula.
This key formula, proved in \cite{HLW}, essentially replaces the
co-area formula applicable in the surface-order asymptotics of
functionals of
sets $A$ having a smooth boundary of bounded curvature \cite{Yu}.

\subsection*{A Steiner-type formula}
Let $A\subset{\Bbb R}^d$ be a non-empty closed set and denote by $\operatorname{exo}(A)$ the exoskeleton of $A$, that is, the set of all $x\in{\Bbb
R}^d\setminus A$ which do not have a unique nearest point in $A$. Then
Theorem~1G in \cite{Fremlin} says that $\mathcal{H}^d(\operatorname{exo}(A))=0$.
Thus, $\mathcal{H}^d$-almost every point $x$ in ${\Bbb R}^d\setminus A$
has a unique nearest point in $A$, denoted by $\pi_A(x)$. The
(reduced) normal bundle $N(A)\subset{\Bbb R}^d\times{\Bbb S}^{d-1}$
of $A$ is given by
\[
N(A):= \biggl\{ \biggl(\pi_A(x),{x-\pi_A(x)\over\llVert  x-\pi_A(x)\rrVert   } \biggr):x
\in{\Bbb R}^d\setminus \bigl(A\cup\operatorname{exo}(A) \bigr) \biggr\},
\]
where here and below $\llVert   \cdot \rrVert   $ stands for the usual Euclidean
distance and $\mathbb{S}^{d-1}$ stands for the Euclidean unit sphere
in $\mathbb{R}^d$. Lemma~2.3 in \cite{HLW} implies that $N(A)$ is a
countably $(d-1)$-rectifiable subset of ${\Bbb R}^d\times{\Bbb
S}^{d-1}$ in the sense of Federer \cite{Federer}, Paragraph 3.2.14.

Let $A$ be as above. The reach function of $A$ is a strictly positive
function on $N(A)$ defined as
\[
\delta(A,x,n):=\inf \bigl\{r\geq0:x+rn\in\operatorname{exo}(A) \bigr\}
\]
for all $(x,n)\in N(A)$. The reach of $A$ is
\[
\operatorname{reach}(A):=\inf \bigl\{\delta(A,x,n):(x,n)\in N(A) \bigr\}
\]
with the convention that $\operatorname{reach}(A)=+\infty$ if $\delta
(A,x,n)=+\infty$ for all $(x,n)\in N(A)$. The set $A$ is said to be of
positive reach if $\operatorname{reach}(A) \in(0, + \infty]$. In
particular, if
$A$ is convex, then \mbox{$\operatorname{reach}(A)=+\infty$}, and vice
versa. We also
remark that any compact $d$-manifold with $C^2$-smooth boundary has
positive reach; cf. \cite{HLW}.

If $A^*$ denotes the closure of the complement of $A$, we see that
$N(\partial A):=N(A)\cup N(A^*)$ and we define the extended normal
bundle of $A$ as $N_e(A):=N(A)\cup TN(A^*)$, where $T$ is the
reflection map $T:\mathbb{R}^d\times{\Bbb S}^{d-1}\rightarrow\mathbb
{R}^d\times
{\Bbb S}^{d-1}, (x,n)\mapsto(x,-n)$. Further, denote the reach
function of $A$ in this context by $\delta^+(A, \cdot, \cdot
)\in[0,+\infty]$ and define the interior reach function $\delta
^-(A,x,n):=-\delta(A^*,x,-n)\in[-\infty,0]$ for $(x,n)\in\mathbb{R}
^d\times{\Bbb S}^{d-1}$.

From Theorem~5.2 in \cite{HLW}, we know that for each $A$ as above
there exist uniquely determined signed measures $\nu_0,\ldots,\nu
_{d-1}$ on $\mathbb{R}^d\times{\Bbb S}^{d-1}$, the so-called extended
support measures of $A$, vanishing outside of $N_e(A)$, such that the
Steiner-type formula
\begin{equation}
\label{eq:Steiner} \int_{\mathbb{R}^d\setminus\partial A}f(x) \,\mathrm{d}x=\sum
_{j=0}^{d-1}\omega _{d-j}\int
_{N_e(A)}\int_{\delta^-(A,x,n)}^{\delta
^+(A,x,n)}r^{d-j-1}
f(x+rn) \,\mathrm{d}r\nu_j \bigl(\mathrm{d}(x,n) \bigr)
\end{equation}
holds for any non-negative measurable bounded function $f:\mathbb
{R}^d\to\mathbb{R}
$ with compact support. Here, for integers $j\geq0$, $\omega_j= j
\kappa
_j:= 2\pi^{j/2}/\Gamma(j/2)$ stands for the %$j$-dimensional
%Hausdorff measure
surface content of the $j$-dimensional unit sphere. The signed measures
$\nu_0,\ldots,\nu_{d-1}$ encode in some sense the singularities of
the boundary of $A$. Although this is not visible in our notation, we
emphasize that the measures $\nu_0,\ldots,\nu_{d-1}$ depend on $A$.

\subsection*{Admissible sets}
Following \cite{HLW}, we denote by
\[
\partial^{+}A:= \bigl\{x\in\partial A:(x,n)\in N(A) \mbox{ for some } n
\in {\Bbb S}^{d-1} \bigr\}
\]
the {positive boundary} of $A$ and define $\operatorname{Nor}(A,x):=\{
n\in{\Bbb
S}^{d-1}:(x,n)\in N(A)\}$ for $x\in\partial^{+} A$. The normal cone
at $x\in\partial^{+}A$ is then $\operatorname{nor}(A,x):=\{an:a\geq0,
n\in\operatorname{Nor}(A,x)\}$ and we put
\begin{equation}
\label{eq:DefA++} \partial^{++}A:= \bigl\{x\in\partial^+A:\operatorname{dim}
\bigl( \operatorname{nor}(A,x) \bigr)=1 \bigr\},
\end{equation}
where $\operatorname{dim}(B)$ denotes the dimension of the affine hull of a set
$B \subset\mathbb{R}^d$.
Clearly, $\partial^{++}A$ is the disjoint union of $\partial^1A$ and
$\partial^2A$, where
\begin{equation}
\label{eq:DefPartialKA} \partial^{k}A= \bigl\{x\in\partial^{++}A:\operatorname{card} \bigl(\operatorname {Nor}(A,x) \bigr)=k \bigr\}, \qquad k\in\{1,2\}.
\end{equation}
Let us recall from \cite{KiderlenRataj} that a closed set $A\subset
\mathbb{R}
^d$ is called {\em gentle} if:
\begin{itemize}[(ii)]
\item[(i)] ${\mathcal{H}}^{d-1}(N_e(A)\cap(B\times\mathbb
{S}^{d-1}))<\infty
$ for all bounded Borel sets $B\subset\mathbb{R}^d$,
\item[(ii)] for ${\mathcal{H}}^{d-1}$-almost all $x\in\partial A$
there are
non-degenerate balls $B_1$ and $B_2$ containing $x$ and satisfying
$B_1\subset A$ and $\operatorname{int}(B_2)\subset\mathbb{R}^d\setminus A$, where
$\operatorname{int}(B_2)$ stands for the interior of $B_2$.
%\item[(i)] $A$ has interior points and the boundary $\partial A$ of
%$A$ is $(\cH^{d-1},d-1)$ rectifiable in the sense of \cite[Paragraph
%3.2.14]{Federer},
%\item[(ii)] $\cH^{d-1}(\partial A\setminus\partial^{++}A)=0$ and $
%\cH^{d-1}(\partial^2A)=0$,
%\item[(iii)] $\left\vert \nu_j\right\vert  (N_e(A))<\infty$ for all $j\in\{0,\ldots,d-1\}$,
\end{itemize}
These assumptions ensure, for example, that ${\mathcal
{H}}^{d-1}(\partial
A\setminus\partial^+A)=0$; cf. equation (5) in \cite{KiderlenRataj}.
%Now, it is known from \cite{HLW} that
The positive boundary of any closed subset of $\mathbb{R}^d$ is
$({\mathcal{H}}
^{d-1},d-1)$-rectifiable \cite{HLW}, and thus the boundary of every
gentle set is $({\mathcal{H}}^{d-1},d-1)$-rectifiable, too. In other words,
there are Lipschitz maps $f_i: \mathbb{R}^{d-1} \to\mathbb{R}^d,  i
= 1,2,\ldots$
such that ${\mathcal{H}}^{d-1}( \partial A \setminus\bigcup_{ i \geq1}
f_i(\mathbb{R}
^{d-1}) ) = 0$; see, for example, \cite{Federer}, Paragraph 3.2.14.
In particular, at ${\mathcal{H}}^{d-1}$-almost every $x\in\partial A$
there is
a unique tangent hyperplane denoted by $T_x:=T_x(\partial A)$.

Moreover, we recall from \cite{KiderlenRataj} that the extended
support measures $\nu_j$ of gentle sets have locally finite total
variation measures $\llvert  \nu_j\rrvert   $ for all $j\in\{0,\ldots,d-1\}$. In
particular, $\llvert  \nu_j\rrvert   (N_e(A))<\infty$ if $A$ is compact.

We now define the class $\mathbf{A}(d)$ of {\em admissible} sets to be the
class of compact sets $A\subset\mathbb{R}^d$ which are gentle,
regular closed
and satisfy ${\mathcal{H}}^{d-1}(\partial^2A)=0$.
(Recall that a set is regular closed if it coincides with the closure
of its interior.)
Here, the assumption that ${\mathcal{H}}^{d-1}(\partial^2A)=0$
simplifies the
structure of the measure $\nu_{d-1}$, to be exploited later.
Regularity excludes sets with lower-dimensional `tentacles' attached
(e.g., a ball with attached line segments). %We notice that each $A\in
%\bA(d)$ has a boundary
%$\partial A$ which is \Comment{? } $(\cH^{d-1},d-1)$ rectifiable in
%the sense that there are Lipschitz maps $f_i: \R^{d-1} \to\R^d,  i =
%1,2,...$ such that
%$\cH^{d-1}( \partial A \setminus\cup_{ i \geq1} f_i(\R^{d-1}) ) =
%0$, see e.g. \cite{Federer}, Paragraph 3.2.14. In particular, at $
%\cH^{d-1}$-almost every $x\in\partial A$ %there is a unique tangent
%hyperplane denoted by $T_x:=T_x(\partial A)$.

The class of gentle and compact sets is rather general and the support
measures $\nu_j$ of such sets
simplify to well-known objects in special situations. We introduce the
following classes of sets:
\begin{itemize}
\item$\mathcal{K}^d$ is the class of convex bodies in $\mathbb
{R}^d$, that is,
compact convex sets $A\subset\mathbb{R}^d$ with non-empty interior,
\item$\mathcal{R}^d$ is the convex ring, consisting of finite unions of
convex bodies in $\mathbb{R}^d$,
\item$\mathcal{M}^d$ denotes the class of compact $d$-dimensional
manifolds in $\mathbb{R}^d$ with twice differentiable boundary,
\item$\mathcal{P}^d$ is the family of compact sets $A\subset\mathbb
{R}^d$ with
positive reach having non-empty interior,
\item$\mathcal{UP}^d$ denotes the class of all subsets $A=A_1\cup
\cdots
\cup A_n\subset\mathbb{R}^d$, $n\geq1$, for sets $A_1,\ldots,A_n\in
\mathcal{P}^d$ and such that $\bigcap_{i\in I}A_i\in\mathcal{P}^d$
for any
$I\subset\{1,\ldots,n\}$.
\end{itemize}
These classes satisfy the inclusions: $\mathcal{K}^d\subset\mathcal{P}^d$,
$\mathcal{K}^d\subset\mathcal{R}^d$, $\mathcal{P}^d\subset\mathcal
{UP}^d$, $\mathcal{M}^d\subset\mathcal{P}^d$ and \mbox{$\mathcal
{R}^d\subset\mathcal{UP}^d$}. If $A\in
\mathcal{K}^d$, then the extended support measures $\nu_j$ are
related to
the generalized curvature measures of $A$ considered in convex
geometry; cf. \cite{SW}. A similar comment applies if $A\in\mathcal
{P}^d$ is a set with positive reach, for which curvature measures have
been introduced in \cite{FedererPaper}. In both cases, it holds that
$\partial^+A=\partial A$. If $A\in\mathcal{K}^d$ then $A$ satisfies
${\mathcal{H}}
^{d-1}(\partial^2A)=0$. The set classes $\mathcal{K}^d$ and $\mathcal{P}^d$
only contain gentle sets. For the set classes $\mathcal{R}^d$ and
$\mathcal{UP}^d$, curvature measures are defined by additive
extension, while for
$\mathcal{M}^d$ curvature measures are defined via classical
differential-geometric methods; see Section~3 in \cite{HLW} for a
detailed discussion. Moreover, for sets $A\in\mathcal{UP}^d$ we have that
${\mathcal{H}}^{d-1}(\partial A\setminus\partial^+A)=0$ (see \cite{HLW}, page~251). Furthermore, if $A\in\mathcal{R}^d$ is regular closed, then
$A$ is gentle according to \cite{KiderlenRataj}, Proposition~2.
Additionally, many $\mathcal{UP}^d$-sets (namely those admitting a
so-called non-osculating representation) are gentle by Proposition~3 in
\cite{KiderlenRataj}.

\subsection*{Admissible score functions}
We next consider the collection $\Xi$ of admissible score functions.
By this we mean the collection of all real-valued Borel measurable
functions $\xi(x, \mathcal{X}, \partial A)$ defined on triples $(x,
\mathcal{X},
\partial A)$, where $\mathcal{X}\subset\mathbb{R}^d$ is locally finite,
$x \in\mathcal{X}$, $A \in\mathbf{A}(d)$, and such that $\xi$ is
translation and
rotation invariant. By the latter two properties, we respectively mean
that $ \xi(x, \mathcal{X}, \partial A) = \xi(x + z, \mathcal{X}+ z,
\partial A+ z)$
and that $ \xi(x, \mathcal{X}, \partial A) = \xi(\vartheta x,
\vartheta\mathcal{X},
\vartheta(\partial A))$ for all $z \in\mathbb{R}^d$, rotations
$\vartheta
\in \operatorname{SO}(d)$ and input $(x, \mathcal{X},\partial A)$. If $x \notin
\mathcal{X}$, we
abbreviate $\xi(x, \mathcal{X}\cup\{x\}, \partial A )$ by
$\xi(x, \mathcal{X}, \partial A )$.

We recall now the concept of a stabilizing functional which was
introduced in \cite{PenroseYukich01,PenroseYukich03,PenroseYukich05} after
earlier works \cite{KestenLee,Lee}; see also the surveys \cite
{Schreiber,YuLMN}. %
%\cite{BaryshnikovYukich2005,Penrose07a,Penrose07b,Schreiber,YuLMN}.. Building on earlier works of Kesten and Lee \cite{KestenLee,Lee},
%Penrose and Yukich introduced stabilizing functionals %
%\cite{PenroseYukich01,PenroseYukich03,PenroseYukich05}, see also
%\cite{BaryshnikovYukich2005,Penrose07a,Penrose07b,Schreiber,YuLMN}.
Roughly speaking, a functional stabilizes if its value at a given point
only depends on a local random neighbourhood and is unaffected by changes
in point configurations outside of it. Following \cite{Yu}, we need to
go a step further in the standard framework
to account for the dependency of functionals $\xi\in\Xi$ on surfaces.
%by taking into account the additional effect that our problems also
%involve a surface $\partial A$.

To make this precise, denote by $B_r(x)$ the Euclidean ball of radius
$r \in(0, \infty)$ and centre $x\in\mathbb{R}^d$ and by $\mathcal
{P}^{\mathrm{hom}}_\tau$ a
homogeneous Poisson point processes on $\mathbb{R}^d$ of intensity
$\tau\in
(0, \infty)$. Say that $\xi\in\Xi$ is homogeneously stabilizing if
for all $\tau\in(0, \infty)$ and all $(d-1)$-dimensional hyperplanes
$H$, there is an almost surely finite random variable $R:=R(\xi,\mathcal{P}^{\mathrm{hom}}
_\tau,H)$ depending on $\xi$, $\mathcal{P}^{\mathrm{hom}}_\tau$
and $H$, the so-called
radius of stabilization, such that
\begin{equation}
\label{eq:hom} \xi \bigl(\mathbf{0}, \mathcal{P}^{\mathrm{hom}}_\tau
\cap B_R(\mathbf {0}), H \bigr) = \xi \bigl(\mathbf{0}, \bigl(
\mathcal{P}^{\mathrm{hom}}_\tau\cap B_R(\mathbf{0})
\bigr) \cup\mathcal{A}, H \bigr)
\end{equation}%
for all locally finite sets $\mathcal{A}\subset B_R(\mathbf{0})^c$,
where $\mathbf{0}$
stands for the origin in $\mathbb{R}^d$. Given \eqref{eq:hom}, the definition
of $\xi$ extends to
Poisson input on all of $\mathbb{R}^d$, that is,
\[
\xi \bigl(\mathbf{0}, \mathcal{P}^{\mathrm{hom}}_\tau, H \bigr) =
\lim_{r \to
\infty} \xi \bigl(\mathbf{0}, \mathcal{P}^{\mathrm{hom}}_\tau
\cap B_r(\mathbf{0}), H \bigr).
\]

Given $A \in\mathbf{A}(d)$, say that $\xi$ is exponentially stabilizing
with respect to the pair $({\mathcal{P}}_\lambda, {\partial A})$ if
for all $x \in\mathbb{R}^d$
there is a random variable $R:=R(\xi,x, {\mathcal{P}}_\lambda,
{\partial A})$, also called
a radius of stabilization, taking values in $[0,\infty)$ with
probability one, such that
\begin{equation}
\label{eq:expo} \xi_\lambda \bigl(x, {\mathcal{P}}_\lambda \cap
B_{\lambda
^{-1/d}R}(x), {\partial A} \bigr) = \xi_\lambda \bigl(x, \bigl({
\mathcal{P}} _\lambda \cap B_{\lambda ^{-1/d}R}(x) \bigr) \cup\mathcal{A}, {
\partial A} \bigr)
\end{equation}
for all locally finite $\mathcal{A}\subset\mathbb{R}^d \setminus
B_{\lambda
^{-1/d}R}(x)$, and the tail probability satisfies
\[
\limsup_{t\to\infty}{1\over t}\log\sup
_{\lambda > 0, x \in
\mathbb{R}^d}\mathbb{P} \bigl[R(\xi,x, {\mathcal{P}}_\lambda,
{\partial A}) > t \bigr]<0.
\]

Surface-order growth for the sums \eqref{eq:gsum} involves finiteness
of the integrated score $\xi_\lambda (x + r\lambda ^{-1/d} n,
{\mathcal{P}}_\lambda, {\partial A})$
over $r \in\mathbb{R}$. Thus, it is natural to require the following
condition; see \cite{Yu}. Given $A \in\mathbf{A}(d)$ and $p \in[1,
\infty
)$, say that $\xi$ satisfies the {\em$p$th moment condition} with
respect to ${\partial A}$ if there is a bounded integrable function
$G^{\xi,p}:= G^{\xi, p, {\partial A}}:
\mathbb{R}\to\mathbb{R}^+$ with $\int_{-\infty}^{\infty} r^{d-1}
(G^{\xi,
p}(r))^{1/p}  \,\mathrm{d}r < \infty$ and such that for all $r \in
\mathbb{R}$ we have
\begin{equation}
\label{eq:mom} \sup_{z \in\mathbb{R}^d \cup\varnothing} \sup_{(x,n)\in N_e(A)} \sup
_{\lambda > 0} \mathbb{E}\bigl\vert \xi_\lambda \bigl(x + r \lambda
^{-1/d} n, {\mathcal{P}}_\lambda \cup z, {\partial A} \bigr)
\bigr\vert ^p \leq G^{\xi,p} \bigl(\llvert r\rrvert \bigr).
\end{equation}

Given $A \in\mathbf{A}(d)$, recall for ${\mathcal{H}}^{d-1}$-almost
all $x \in{\partial A}$ that
$T_x:=T_x(\partial A)$ is the unique hyperplane tangent to ${\partial
A}$ at $x$.
%Put $H_x:=H(\0,\pA-x)$.
For $x \in{\partial A}$, we put $H_x:= T_{\mathbf{0}}({\partial A}- x)$.
The score $\xi$ is said to be {\em well approximated by ${\mathcal
{P}}_\lambda $
input on half-spaces} if for all $A \in\mathbf{A}(d)$, almost all $x
\in{\partial A}
$, and all $w \in\mathbb{R}^d$, we have
\begin{equation}
\label{eq:lin} \lim_{\lambda \to\infty} \mathbb{E}\bigl\vert \xi \bigl(w,
\lambda^{1/d}({\mathcal{P}}_\lambda - x),\lambda
^{1/d}({\partial A}- x) \bigr) - \xi \bigl(w, \lambda ^{1/d}({
\mathcal{P}}_\lambda - x), H_x \bigr)\bigr\vert = 0.
\end{equation}

\subsection*{General theorems giving first- and second-order asymptotics}
The results asserted in Section~\ref{sect1} are consequences of general limit
theorems giving expectation
and variance asymptotics for the statistics \eqref{eq:gsum}. We first
describe the general theory and
then, in Section~\ref{sec:ProofsApplications}, show how to deduce the assertions of Section~\ref{sect1}.
The general limit theorems given here extend Theorems 1.1 and 1.2 in
\cite{Yu} to the class of admissible sets
and they yield the first- and second-order asymptotics for statistics
of other surfaces, as discussed in Remark (iii) below.

%We now state two general limit theorems concerning expectation and
%variance asymptotics for the statistics of the form \eqref{eq:gsum}.
%The results asserted in Section~1 are then deduced from our general
%limit theorems, as shown in Section~4.
%\Comment{ } The general limit theorems given here extend Theorems 1.1
%and 1.2 in \cite{Yu} to the class of admissible sets.
% We anticipate that the general results yield the limit theory for
%statistics of other surfaces as discussed in Remark (iii) after
%Theorem~\ref{eq:Var} below.

For a score function $\xi\in\Xi$, we put
\begin{equation}
\label{expec-formula} \mu(\xi,{\partial A}):=\int_{\partial^1 A}\int
_{-\infty}^\infty {\Bbb E}\xi \bigl(\mathbf{0}+sn,
\mathcal{P}_{\kappa(x)}^\mathrm{hom}, \mathbb {R}^{d-1} \bigr)
\kappa (x) \,\mathrm{d}s \mathcal{H}^{d-1}(\mathrm{d}x),
\end{equation}
where $n$ is the unique unit normal at $\mathbf{0}$ with respect to
$\mathbb{R}^{d-1}$.
%we recall from \eqref{eq:DefPartialKA} the definition of $\partial^1
%A$ that $n(x)$ is the unique unit normal at $x\in\partial^1 A$.,
%recall \eqref{eq:DefPartialKA}.
%$\cH^{d-1}$-almost everywhere unique unit normal at $x\in\pA$, recall
%\eqref{eq:DefPartialKA}. %Also recall our standing assumption that $
%\xi$ is translation and rotation invariant. %We are now prepared to
%state our first general result, which provides expectation asymptotics.
We now state a general result giving expectation asymptotics for sums
of score functions.
Let $\mathcal{C}({\partial A})$ denote the set of functions on
$\mathbb{R}^d$ which are
continuous at all points $x \in{\partial A}$.

\begin{theorem} \label{eq:WLLN}
Let $A \in\mathbf{A}(d)$ and $\kappa \in\mathcal{C}({\partial
A})$. Suppose that $\xi\in
\Xi$ is homogeneously stabilizing \eqref{eq:hom}, satisfies the
moment condition \eqref{eq:mom} for some $p \in[1, \infty)$, and is
well approximated by ${\mathcal{P}}_\lambda $ input on half-spaces as
at \eqref
{eq:lin}. Then for $m \in\{1,2\}$, we have the following weak law of
large numbers:
\begin{equation}
\label{expec}
 \lim_{\lambda\to\infty}
\lambda^{-(d-1)/d} H^\xi(\mathcal {P}_\lambda,\partial A) =
\mu(\xi,{\partial A})\qquad {\mbox{in }} L^m.
\end{equation}
\end{theorem}

%\CT[Remark on $L^2$ convergence. Shall we really give the details?]
%\begin{remark}
%Following the proof of Theorem~1.1 in \cite{Yu}, one can extend
%Theorem~\ref{eq:WLLN} from $L^1$- to $L^2$-convergence.
%\end{remark}

Next, we turn to variance asymptotics and define for $x, x' \in\mathbb{R}^d$,
$\tau\in(0, \infty)$, and all $(d-1)$-dimensional hyperplanes $H$,
\begin{eqnarray*}
c^\xi \bigl(x,x'; \mathcal{P}^{\mathrm{hom}}_\tau,
H \bigr) &:= & \mathbb{E}\xi \bigl(x, \mathcal{P}^{\mathrm{hom}}_\tau
\cup \bigl\{x' \bigr\}, H \bigr) \xi \bigl(x',
\mathcal{P}^{\mathrm{hom}}_\tau\cup\{x\}, H \bigr)
\\
&&{} - \mathbb{E}\xi \bigl(x, \mathcal{P}^{\mathrm{hom}} _\tau, H \bigr)
\mathbb{E}\xi \bigl(x', \mathcal{P}^{\mathrm{hom}}_\tau,
H \bigr).
\end{eqnarray*}
Moreover, define $\sigma^2(\xi, {\partial A})$ by
\begin{eqnarray}
\label{eq: sigma}
  &&\sigma^2(
\xi, {\partial A})\nonumber
\\
 &&\quad:= \mu \bigl(\xi^2, {\partial A} \bigr)
\\
&&\qquad {}+ \int_{\partial^1A} \int_{\mathbb{R}^{d-1}} \int
_{-\infty
}^{\infty}\int_{-\infty}^{\infty}
c^\xi \bigl(\mathbf{0} + rn, p+sn; {\mathcal{P}}_{\kappa (x)}^\mathrm{hom},
\mathbb{R}^{d-1} \bigr) \kappa (x)^2 \,\mathrm{d}s
\,\mathrm{d}r \,\mathrm{d}p \mathcal{H}^{d-1}(\mathrm{d}x).\qquad \nonumber
\end{eqnarray}
The following general result gives variance asymptotics for sums of
score functions.

\begin{theorem} \label{eq:Var}
Let $A \in\mathbf{A}(d)$ and $\kappa \in\mathcal{C}({\partial
A})$. We assume that $\xi
\in\Xi$ is homogeneously stabilizing \eqref{eq:hom},
exponentially stabilizing \eqref{eq:expo}, satisfies the moment
condition \eqref{eq:mom} for some $p \in(2, \infty)$
%with $G^{\xi,2}$ decaying exponentially fast,
and is well approximated by ${\mathcal{P}}_\lambda $ input on
half-spaces as at
\eqref{eq:lin}. Then
\begin{equation}
\label{eq:Var1} \lim_{\lambda \to\infty} \lambda ^{-(d-1)/d}
\operatorname {Var} \bigl[H^\xi({\mathcal{P}}_\lambda, {
\partial A} ) \bigr] = \sigma^2(\xi, {\partial A}).
\end{equation}
\end{theorem}

Some of the applications presented in Section~\ref{sect1} require the limit
theory for the non-re-scaled sums $\sum_{x \in{\mathcal{P}}_\lambda
} \xi(x, {\mathcal{P}}
_\lambda, \partial A)$. To state the result in this case,
%ecall \cite{YuLMN}
call a score function $\xi$ {\em homogeneous of order $\gamma\in
\mathbb{R}
$} if for all $a \in(0, \infty)$,
\[
\xi \bigl(ax, a\mathcal{X}, a({\partial A}) \bigr)= a^{\gamma} \xi(x,
\mathcal{X}, {\partial A}).
\]
When $\xi$ is homogeneous of order $\gamma$, it follows that
\[
\sum_{x \in{\mathcal{P}}_\lambda } \xi(x, {\mathcal{P}}_\lambda,
\partial A) = \lambda ^{ {-\gamma
/ d}} H^\xi({\mathcal{P}}_\lambda, {\partial A}).
\]
Homogeneity, together with the distributional identity ${\mathcal
{P}}_{\kappa (x)
}\stackrel{\mathcal{D}}{=}\kappa (x)^{-1/d}  {\mathcal{P}}_1 $
gives
\begin{eqnarray}
\nonumber
\mu(\xi,{\partial A}) &=& \int_{-\infty}^\infty{
\Bbb E}\xi \bigl(\mathbf{0}+sn,\mathcal{P}_{1}^\mathrm{hom},
\mathbb{R}^{d-1} \bigr) \,\mathrm{d}s \int_{\partial^1 A} \kappa
(x)^{1 - \gamma/d} \mathcal {H}^{d-1}(\mathrm{d} x)\nonumber\label{kad}
\\[-8pt]\\[-8pt]
&=& \int_{-\infty}^\infty{\Bbb E}\xi \bigl(\mathbf{0}+sn,
\mathcal {P}_{1}^\mathrm{hom}, \mathbb{R}^{d-1} \bigr) \,\mathrm{d}s \cdot{\mathcal{H}}_{\kappa,\gamma}^{d-1}(\partial A)\nonumber
\end{eqnarray}
and
\begin{eqnarray}
 &&\sigma^2(\xi, {\partial A})\nonumber
 \\
  &&\quad = \int
_{-\infty}^\infty{\Bbb E} \xi ^2 \bigl(\mathbf{0}+sn, \mathcal{P}_{1}^\mathrm{hom}, \mathbb{R}^{d-1}
\bigr) \,\mathrm{d}s \int_{\partial^1 A} \kappa (x)^{1 - 2\gamma/d} \mathcal
{H}^{d-1}(\mathrm{d} x)\nonumber
\\
&&\qquad {} + \int_{\mathbb{R}^{d-1}} \int_{-\infty}^{\infty}
\int_{-\infty
}^{\infty} c^{\xi} \bigl(\mathbf{0} +
rn, p+sn; {\mathcal{P}}_{1}^\mathrm{hom}, \mathbb{R}^{d-1}
\bigr) \,\mathrm{d}s \,\mathrm{d}r \,\mathrm{d}p\nonumber
\\[-8pt]\\[-8pt]
&&\qquad {}\times  \int_{\partial^1A} \kappa
(x)^{2 -
2\gamma
/d} \mathcal{H}^{d-1}(\mathrm{d}x)\nonumber
\\
\label{vkad}
& &\quad  =\int_{-\infty}^\infty{
\Bbb E} \xi^2 \bigl(\mathbf{0}+sn, \mathcal {P}_{1}^\mathrm{hom}, \mathbb{R}^{d-1} \bigr) \,\mathrm{d}s \cdot{\mathcal{H}}_{\kappa,2\gamma}^{d-1}(
\partial A)\nonumber
\\
&&\qquad {} + \int_{\mathbb{R}^{d-1}} \int_{-\infty}^{\infty}
\int_{-\infty
}^{\infty} c^{\xi} \bigl(\mathbf{0} +
rn, p+sn; {\mathcal{P}}_{1}^\mathrm{hom}, \mathbb{R}^{d-1}
\bigr) \,\mathrm{d}s \,\mathrm{d}r \,\mathrm{d}p \cdot{\mathcal{H}}_{\kappa ^2,\gamma
}^{d-1}(
\partial A).\nonumber
\end{eqnarray}
Consequently, with $\mu(\xi,{\partial A})$ and $\sigma^2(\xi,
{\partial A})$ as
in \eqref{kad} and \eqref{vkad}, respectively, we have under the
conditions of Theorems \ref{eq:WLLN} and \ref{eq:Var} that
\begin{equation}
\label{eq: newexpect} \lim_{\lambda \to\infty} \lambda ^{-(d-1 - \gamma)/d}\sum
_{x \in
{\mathcal{P}}_\lambda } \xi(x, {\mathcal{P}}_\lambda,
\partial A) = \mu(\xi,{\partial A})
\end{equation}
in $L^m$ for $m\in\{1,2\}$, and
\begin{equation}
\label{eq: newvar} \lim_{\lambda \to\infty} \lambda ^{-(d-1 - 2\gamma)/d}
\operatorname{Var}\sum_{x \in
{\mathcal{P}}_\lambda } \xi(x, {
\mathcal{P}}_\lambda, \partial A) = \sigma^2(\xi, {\partial
A}).
\end{equation}

\begin{remarks*}
\begin{longlist}
\item[(i)] {\em Convergence of random measures.} The methods presented
here also yield expectation and variance asymptotics for integrals of
the empirical measures
\[
\sum_{x \in{\mathcal{P}}_\lambda } \xi_\lambda (x, {\mathcal
{P}}_\lambda, \partial A) \delta_x
\]
against elements of $\mathcal{C}({\partial A})$ (here, $\delta_x$
stands for the
unit-mass Dirac measure at $x$).
The details of this extension are straightforward and may be found in,
for example, \cite{YuLMN}, which deals with volume-order asymptotics
for sums of score functions.

\item[(ii)] {\em Central limit theorems.} Say that $\xi$ decays {\em
exponentially fast} with respect to the distance to $\partial A$ if for
all $p \in[1, \infty)$ the function $G^{\xi,p}$ defined at \eqref
{eq:mom} satisfies
\begin{equation}
\label{strongmom} \limsup_{\llvert  u\rrvert    \to\infty}\llvert u\rrvert ^{-1}
\log G^{\xi,p} \bigl(\llvert u\rrvert \bigr) < 0.
\end{equation}
Let $\Phi( \cdot )$ denote the distribution
function of a standard normal random variable.
If $\xi\in\Xi$ decays exponentially fast as in \eqref{strongmom}
and if $\xi$ satisfies the moment condition \eqref{eq:mom} with
$p=3$, then
by Theorem~1.3 of \cite{Yu}, the statistics \eqref{eq:gsum} satisfy a
central limit theorem
\[
\sup_{x\in\mathbb{R}} \biggl\llvert P \biggl[{H^\xi({\mathcal{P}}_\lambda,{\partial A})-\mathbb{E}H^\xi({\mathcal{P}}_\lambda,{\partial A})\over\sqrt{\operatorname{Var}[H^\xi({\mathcal
{P}}_\lambda,{\partial A})]}}
\biggr]-\Phi(x) \biggr\rrvert \leq r(\lambda )
\]
with rate function
\[
r(\lambda ):=c(\log\lambda )^{3d+1}\lambda ^{(d-1)/d} \bigl(
\operatorname {Var} \bigl[H^\xi(P\lambda,{\partial A} ) \bigr]
\bigr)^{3/2},
\]
where $c>0$ is a constant not depending on $\lambda $. In particular, if
$\sigma^2(\xi,\partial A)$ is strictly positive, then $r(\lambda ) =
c(\log\lambda )^{3d+1}\lambda ^{-(d-1)/2d}$. This is the case for
the examples in Section~\ref{sect1}, provided that $\kappa\equiv1$ and
that $\partial A$ contains a $C^2$-smooth subset with positive
${\mathcal{H}}
^{d-1}$-measure.

%By Theorem~1.3 of \cite{Yu}, the statistics \eqref{eq:gsum} satisfy a
%central limit theorem. If $\Phi( \cdot )$ denotes the distribution
%function of a standard normal %random variable and if $\xi\in\Xi$
%satisfies the moment condition \eqref{eq:mom} with $p=3$, then $$
%\sup_{x\in\R}\Big\vert P\Big[{H^\xi(\P_\la,\pA)-\E H^\xi(\P_\la,\pA)\over%
%\sqrt{\Var[H^\xi(\P_\la,\pA)]}}\Big]-\Phi(x)\Big\vert \leq r(\la)$$ with
%rate function $$r(\la):=C(\log\la)^{3d+1}\la^{(d-1)/d}(\Var[H^\xi(P\la,
%\pA)])^{3/2},$$ where $C>0$ is a %constant not depending on $\la$.
%In particular, if $\sigma^2(\xi,\partial A)$ is strictly positive,
%then $r(\la) = C((\log\la)^{3d+1}\la^{-(d-1)/2d})$. This is the case
%for %the examples presented in the introduction, provided that $
%\partial A$ contains a $C^2$-smooth subset with positive $
%\cH^{d-1}$-measure and if $\kappa\equiv1$.

\item[(iii)] {\em Further applications of general results.} Theorems
\ref{eq:WLLN} and \ref{eq:Var} have applications
to statistics of surfaces going beyond those arising in
Poisson--Voronoi approximation. For instance, these general theorems
provide the limit theory for functionals of surfaces of germ--grain
models including, for example, the limit theory for the number of
exposed tangent points to Boolean models, as described in Section~3.2
of \cite{Mo}. Another application of the general theory
involves the number of {\em maximal points} in a sample, which goes as follows.
% We indicate the problem of so-called {\em maximal points}.
A point $x\in{\mathcal{P}}_{\lambda }$ is called maximal if the
Minkowski sum $(\mathbb{R}
_+)^d\oplus x$ contains no other point of ${\mathcal{P}}_{\lambda }$
besides $x$,
that is, if $((\mathbb{R}_+)^d\oplus x)\cap{\mathcal{P}}_{\lambda
}=\{x\}$. The number
$M_\lambda $ of maximal points of ${\mathcal{P}}_{\lambda }$ has
attracted considerable
interest in the literature; see \cite
{BaiHwangEtal01,BarbourXia01,BarbourXia06,Devroye,HwangTsai10,Yu}.
These works restrict to domains $A$ that are either piecewise linear,
convex or smooth. We may use Theorems \ref{eq:WLLN} and \ref{eq:Var}
to unify and extend these results to domains $A$ which are admissible
sets, as illustrated by the following statement, whose proof follows
from modifications of the proof of Theorem~2.5 in \cite{Yu} and is
left to the reader. Let $\kappa$ be a density supported on $A:=\{
(u,v)\in\mathbb{R}^{d-1}\times\mathbb{R}:u\in D, 0\leq v\leq f(u)\}
$, where
$D\subset\mathbb{R}^{d-1}$ and $f:D\to\mathbb{R}$, and assume that
$A$ is an
admissible set, that is, $A \in\mathbf{A}(d)$. We further assume that the
partial derivatives of $f$ exist a.e. and are bounded away from zero
and infinity.
If ${\mathcal{P}}_{\lambda }$ is a Poisson point process whose
intensity measure has
density $\lambda \kappa$ with respect to Lebesgue measure
then there are constants $c_{11} \in(0, \infty)$ and $c_{12} \in[0,
\infty)$ depending only on $d$, $\kappa$ and $A$ such that
\[
\mathbb{E}M_\lambda \sim c_{11}\lambda ^{1-{1/d}} \quad \mbox{and}\quad  \operatorname{Var} [M_\lambda ]\sim c_{12}\lambda
^{1-{1/d}}.
\]
\end{longlist}
\end{remarks*}

%s3
\section{Proofs of Theorems \texorpdfstring{\protect\ref{eq:WLLN}}{2.1} and \texorpdfstring{\protect\ref{eq:Var}}{2.2}}\label
{sec:ProofsGeneralTheorems}

%\subsection{Preparatory lemmas}

To keep the paper self-contained, we give three preparatory lemmas
pertaining to the re-scaled scores $\xi_\lambda, \lambda > 0$.
These are
re-formulations of Lemmas 3.1--3.3 in \cite{Yu}, which we adopt to our
more general set-up. The following lemmas do not require continuity of
$\kappa $ but instead use that a.e. $x \in\mathbb{R}^d$ is a
Lebesgue point of
$\kappa $, that is to say
\[
{1\over\varepsilon^d} \int_{B_{\varepsilon}(x) } \bigl\vert \kappa (y) -
\kappa(x)\bigr\vert \,\mathrm{d}y
\]
tends to zero as $\varepsilon\downarrow0$.
Given $x\in\partial^1 A$, with $\partial^1 A$ defined at \eqref
{eq:DefPartialKA},
recall that $H_x:=T_\mathbf{0}(\partial A-x)$ is the unique tangent
hyperplane to ${\partial A}- x$ at $\mathbf{0}$ with unit normal
$n(x)$. %If $A\in
%\bA(d)$
Let $\mathbf{0}_x$ denote a point at the origin of $H_x$.

\begin{lemma} \label{eq:L1}
Fix $A \in\mathbf{A}(d)$. Assume that $\xi$ is homogeneously
stabilizing as
at \eqref{eq:hom}, satisfies the
moment condition \eqref{eq:mom} for some $p \in(1, \infty)$ and is
well approximated by ${\mathcal{P}}_\lambda $ input on half-spaces
\eqref{eq:lin}.
%%Let $x\in\partial^1 A$ and denote by %$T_x$ the unique tangent
%hyperplane at $x$ with unit normal $n(x)$.
Then for all $x\in\partial^1A$, $w \in\mathbb{R}^d$, and $r\in
\mathbb{R}$ we have
\begin{equation}
\label{eq:L3.1f} \lim_{\lambda \to\infty} \mathbb{E}\xi_\lambda
\bigl(x + r\lambda ^{-1/d}n(x) + \lambda ^{-1/d}w, {
\mathcal{P}}_\lambda, {\partial A} \bigr) = \mathbb{E}\xi \bigl(
\mathbf{0}_x + rn(x) + w, {\mathcal{P}}_{\kappa
(x)}^\mathrm{hom},
H_x \bigr).
\end{equation}
\end{lemma}

\begin{lemma} \label{L2}
Fix $A \in\mathbf{A}(d)$. Assume that $\xi$ is homogeneously
stabilizing as
at \eqref{eq:hom}, satisfies the
moment condition \eqref{eq:mom} for some $p \in(2, \infty)$, and is
well approximated by
${\mathcal{P}}_\lambda $ input on half-spaces \eqref{eq:lin}.
Given $x\in\partial^1 A$, $v \in\mathbb{R}^d$, and $r \in\mathbb
{R}$, %let $T_x$
%be the unique tangent hyperplane at $x$ with unit normal $n(x)$, and
we put for $\lambda \in(0,\infty)$,
\begin{eqnarray*}
X_\lambda &:= &\xi_\lambda \bigl(x + r\lambda
^{-1/d}n(x), {\mathcal {P}}_\lambda \cup \bigl\{x + r\lambda
^{-1/d}n(x) +\lambda ^{-1/d}v \bigr\}, {\partial A} \bigr),
\\
Y_\lambda &:=& \xi_\lambda \bigl(x + r\lambda ^{-1/d}n(x)
+ \lambda ^{-1/d}v, {\mathcal{P}}_\lambda \cup \bigl\{ x + r
\lambda ^{-1/d}n(x) \bigr\}, {\partial A} \bigr),
\\
X &:=& \xi \bigl(\mathbf{0}_x + rn(x), {\mathcal{P}}_{\kappa (x)}^\mathrm{hom}
\cup \bigl\{\mathbf{0}_x + rn(x) + v \bigr\}, H_x
\bigr),
\\
Y &:=& \xi \bigl(\mathbf{0}_x + rn(x) + v, {\mathcal{P}}_{\kappa (x)}^\mathrm{hom}
\cup \bigl\{\mathbf{0}_x + rn(x) \bigr\}, H_x \bigr).
\end{eqnarray*}
Then $ \lim_{\lambda \to\infty} \mathbb{E}[ X_\lambda
Y_\lambda ]= \mathbb{E}[ X Y] $.
\end{lemma}

%\CT[I wrote out the norm, because we never introduce the notation. ok
%JY]

\begin{lemma} \label{eq:L3}
Fix $A \in\mathbf{A}(d)$. Let $\xi$ be exponentially stabilizing as at
\eqref{eq:expo} and assume the moment condition \eqref{eq:mom} holds
for some $p \in(2, \infty)$. Then there is a constant $C\in(0,
\infty)$ such that for all $w, v \in\mathbb{R}^d$ and $\lambda \in
(0, \infty
)$, we have
\begin{eqnarray*}
&& \bigl\vert \mathbb{E}\xi_\lambda \bigl(w, {\mathcal{P}}_\lambda
\cup \bigl\{w + \lambda ^{-1/d}v \bigr\}, {\partial A} \bigr) \xi
_\lambda \bigl(w + \lambda ^{-1/d}v, {\mathcal{P}}_\lambda
\cup\{w\}, {\partial A} \bigr)
\\
&&\qquad {} - \mathbb{E}\xi_\lambda (w, {\mathcal {P}}_\lambda, {\partial
A}) \mathbb{E}\xi _\lambda \bigl(w + \lambda ^{-1/d}v, {
\mathcal{P}}_\lambda, {\partial A} \bigr) \bigr\vert
\\
&&\quad  \leq C \bigl(\mathbb{E}\xi_\lambda \bigl(w, {\mathcal{P}}_\lambda
\cup \bigl\{w + \lambda ^{-1/d}v \bigr\}, {\partial A}
\bigr)^p \bigr)^{1/p}
\\
&&\qquad {} \times \bigl(\mathbb{E}\xi_\lambda \bigl(w + \lambda ^{-1/d}v,
{\mathcal{P}} _\lambda \cup\{w\}, {\partial A} \bigr)^p
\bigr)^{1/p} \exp \bigl(-C^{-1}\llVert v\rrVert \bigr).
\end{eqnarray*}
In particular, there is a constant ${c}\in(0,\infty)$ such that
if $w = x + r\lambda ^{-1/d}n(x)$, then
\begin{eqnarray*}
&& \bigl\vert \mathbb{E}\xi_\lambda \bigl(w, {\mathcal{P}}_\lambda
\cup \bigl\{w + \lambda ^{-1/d}v \bigr\}, {\partial A} \bigr) \xi
_\lambda \bigl(w + \lambda ^{-1/d}v, {\mathcal{P}}_\lambda
\cup\{w\}, {\partial A} \bigr)
\\
&&\quad {} - \mathbb{E}\xi_\lambda (w, {\mathcal{P}}_\lambda, {\partial
A}) \mathbb{E}\xi_\lambda \bigl(w + \lambda ^{-1/d}v, {
\mathcal{P}}_\lambda, {\partial A} \bigr) \bigr\vert \leq{c} G^{\xi, p}
\bigl(\llvert r\rrvert \bigr)^{1/p} \exp \bigl(-{c}^{-1}\llVert
v\rrVert \bigr).
\end{eqnarray*}
%
%with a constant $\bar{c}\in(0,\infty)$.
\end{lemma}

\begin{pf}
The first asserted inequality follows as in either Lemma~4.2 of \cite
{Penrose07a} or Lemma~4.1 of \cite{BaryshnikovYukich2005}.
The second assertion follows from the first assertion together with the
moment condition \eqref{eq:mom}.
\end{pf}

%\subsection{Proofs of the general theorems}

Given these auxiliary lemmas, we may now prove the general results.

\begin{pf*}{Proof of Theorem~\ref{eq:WLLN}}
To show \eqref{expec}, it is enough to show the expectation asymptotics
\begin{equation}
\label{expec1} \lim_{\lambda\to\infty} \lambda^{-(d-1)/d} \mathbb{E}
\sum_{x\in
\mathcal{P}_\lambda}\xi_\lambda(x,\mathcal{P}_\lambda,
\partial A) = \mu(\xi,{\partial A})
\end{equation}
and then follow the method of proof of Theorem~1.1 of \cite{Yu} to
deduce $L^m$-convergence for $m\in\{1,2\}$.

To show \eqref{expec1}, we first apply the Mecke identity \cite{SW}, Theorem~3.2.5, for Poisson point processes to obtain
\begin{eqnarray*}
  \lambda^{-(d-1)/d} {\Bbb E}\sum
_{x\in\mathcal{P}_\lambda}\xi _\lambda(x,\mathcal{P}_\lambda,
\partial A) &=& \lambda ^{-(d-1)/d}\int_{\mathbb{R}^d}{\Bbb E}
\xi_\lambda(x,\mathcal{P}_\lambda,\partial A) \lambda \kappa(x)
\,\mathrm{d}x
\\
&=&\lambda^{1/d}\int_{\mathbb{R}^d}{\Bbb E}
\xi_\lambda(x,\mathcal {P}_\lambda,\partial A) \kappa(x)
\,\mathrm{d}x;   %
\end{eqnarray*}
recall that we write $\xi_\lambda(x,\mathcal{P}_\lambda,\partial A)$
instead of $\xi_\lambda(x,\mathcal{P}_\lambda\cup\{x\},\partial
A)$ if
$x\notin\mathcal{P}_\lambda$.
We now use the Steiner-type formula \eqref{eq:Steiner} to re-write the
last integral as
\[
\lambda^{1/d}\sum_{j=0}^{d-1}
\omega_{d-j}\int_{N_e(A)}\int_{T(x,n)}r^{d-j-1}
{\Bbb E}\xi_\lambda(x+rn,\mathcal{P}_\lambda,\partial A)
\kappa(x+rn) \,\mathrm{d}r \nu_j \bigl(\mathrm{d}(x,n) \bigr),
\]
where for fixed $(x,n)\in N_e(A)$, $T(x,n):= [\delta^-(A,x,n),\delta
^+(A,x,n)]$.
%such that $x+rn\in Q$ (recall once more the r\^ole of the exoskeleton
%of $A$ in the definition of the (extended) normal bundle).
Upon the substitution \mbox{$r=\lambda^{-1/d}r'$}, we obtain that $\lambda
^{-(d-1)/d}  {\Bbb E}\sum_{x\in\mathcal{P}_\lambda}\xi_\lambda
(x,\mathcal{P}_\lambda,\partial A) $ equals
\begin{eqnarray}
\label{INTA}
 && \sum_{j=0}^{d-1}
\omega_{d-j}\lambda^{-(d-1-j)/d}\int_{N_e(A)}\int
_{\lambda^{1/d}T(x,n)}  \bigl(r' \bigr)^{d-j-1} {\Bbb E}
\xi_\lambda \bigl(x+\lambda ^{-1/d}r'n,
\mathcal{P}_\lambda,\partial A \bigr)\nonumber
\\[-8pt]\\[-8pt]
%&\one(s \leq\la^{1/d}T(x,n))
&&\quad {} \times\kappa \bigl(x+\lambda^{-1/d}r'n \bigr)
\,\mathrm{d}r' \nu_j \bigl(\mathrm {d}(x,n) \bigr).\nonumber
\end{eqnarray}
To simplify the notation, write $r$ for $r'$.
By the moment assumption \eqref{eq:mom} with $p = 1$, we conclude
that, for each $j\in\{0,\ldots,d-1\}$, the integrand is bounded by
the product $\llvert  r\rrvert   ^{d-j-1} G^{\xi, 1}(\llvert  r\rrvert   )\llVert  \kappa\rrVert   _\infty$,
implying that
\begin{eqnarray*}
  && \biggl\llvert \int_{N_e(A)} \int
_{\lambda^{1/d}T(x,n)}r^{d-j-1} {\Bbb E}\xi_\lambda \bigl(x+
\lambda^{-1/d}rn,\mathcal{P}_\lambda,\partial A \bigr) \kappa
\bigl(x+\lambda^{-1/d}rn \bigr) \,\mathrm{d}r \nu_j \bigl(
\mathrm{d}(x,n) \bigr) \biggr\rrvert
\\
&&\quad  \leq\int_{N_e(A)}\int_{-\infty}^\infty
r^{d-j-1} G^{\xi,
1} \bigl(\llvert r\rrvert \bigr)\llVert \kappa
\rrVert _\infty \,\mathrm{d}r\llvert \nu_j\rrvert \bigl(
\mathrm{d}(x,n) \bigr)
\\
&&\quad  =\llVert \kappa\rrVert _\infty\llvert \nu_j\rrvert
\bigl(N_e(A) \bigr) \int_{-\infty
}^{\infty}
r^{d-j-1} G^{\xi, 1} \bigl(\llvert r\rrvert \bigr) \,\mathrm{d}r.
\end{eqnarray*}
The integral $\int_{-\infty}^{\infty} r^{d-j-1} G^{\xi, 1}(\llvert  r\rrvert   )
\,\mathrm{d}r$ is finite by assumption. Moreover, $\llVert  \kappa\rrVert   _\infty
<\infty
$ by assumption and $\llvert  \nu_j\rrvert   (N_e(A))<\infty$ since $A\in\mathbf{A}(d)$.
Consequently, taking the limit in \eqref{INTA} as $\lambda \to\infty$,
it follows by the dominated convergence theorem that
only the term $j=d-1$ remains:
\begin{eqnarray}
\label{eq:lim2}
  && \lim_{\lambda\to\infty}
\lambda^{-(d-1)/d} {\Bbb E}\sum_{x\in
\mathcal{P}_\lambda}
\xi_\lambda(x,\mathcal{P}_\lambda,\partial A)\nonumber
\\
&&\quad  = 2\int_{N_e(A)}\int_{-\infty}^\infty
\lim_{\lambda\to\infty
}{\Bbb E}\xi_\lambda \bigl(x+
\lambda^{-1/d}rn,\mathcal{P}_\lambda,\partial A \bigr)
\\
&&\qquad {} \times\kappa \bigl(x+\lambda^{-1/d}rn \bigr)\mathbf{1} \bigl(r \in
\lambda^{1/d}T(x,n) \bigr) \,\mathrm{d}r \nu_{d-1} \bigl(
\mathrm{d}(x,n) \bigr).  \nonumber %
\end{eqnarray}
Here, we use the identity $\omega_1 = 2$ and we also use that $\lim_{\lambda\to\infty} \lambda^{1/d}T(x,n) = (-\infty,\infty
)$, which holds by construction of $N_e(A)$, where the exoskeleton has
been excluded. By continuity of $\kappa $ on~$\partial A$, we have
$\lim_{\lambda\to\infty} \kappa(x+\lambda^{-1/d}rn) = \kappa
(x)$. Finally, consider the limit
\[
\lim_{\lambda\to\infty}{\Bbb E}\xi_{\lambda} \bigl(x+\lambda
^{-1/d}rn,\mathcal{P}_\lambda,\partial A \bigr).
\]
To identify it, we use translation invariance and the definition of
$\xi_\lambda$, and write
\begin{eqnarray*}
 && \xi_{\lambda} \bigl(x+\lambda^{-1/d}rn,
\mathcal{P}_\lambda,\partial A \bigr)
\\
 &&\quad = \xi_\lambda \bigl(\mathbf{0}_x+\lambda^{-1/d}rn,\mathcal{P}_\lambda -x,
\partial A-x \bigr)
\\
&&\quad =\xi \bigl(\mathbf{0}_x+rn,\lambda^{1/d}(
\mathcal{P}_\lambda-x),\lambda ^{1/d}(\partial A-x) \bigr).
\end{eqnarray*}
The measure $\nu_{d-1}$ concentrates, according to the discussion around
Proposition~4.1 of \cite{HLW}, Section~4, on the subset ${\partial
A}^{++}$ of
the boundary $\partial A$ where the normal cone is
one dimensional; recall \eqref{eq:DefA++}. Moreover, since $A\in
\mathbf{A}
(d)$, the measure $\nu_{d-1}$ in fact concentrates on the subset
$\partial^1A\subset\partial^{++}A$ (see \eqref{eq:DefPartialKA}),
that is to say, on points of the boundary having a unique normal vector
or tangent hyperplane as in the case of a smooth surface.

Since $\xi$ is well approximated by input on half-spaces, Lemma~\ref
{eq:L1} implies for all $(x,n) \in N_e(A)$ with $x\in\partial^1A$,
that the expectation of the latter expression converges to
\[
\lim_{\lambda\to\infty} \mathbb{E}\xi \bigl(\mathbf{0}_x+rn,
\lambda ^{1/d}(\mathcal{P}_\lambda-x),\lambda ^{1/d}(
\partial A-x) \bigr) = {\Bbb E}\xi \bigl(\mathbf{0}_x+rn,
\mathcal{P}_{\kappa(x)}^\mathrm{hom}, H_x \bigr).
\]
%
%$${\Bbb E}\xi(\mathbf{0}_x+sn,\mathcal{P}_{\kappa(x)}^\mathrm{hom},{ Tan}(
%\partial%A,x)).$$
Thus, we obtain from \eqref{eq:lim2},
\begin{eqnarray}
\label{eq:AlmostFinishedExpectationAsymptotic}
  & &\lim_{\lambda\to\infty}{
\lambda^{-(d-1)/d}} {\Bbb E}\sum_{x\in
\mathcal{P}_\lambda}
\xi_\lambda(x,\mathcal{P}_\lambda,\partial A)\nonumber
\\[-8pt]\\[-8pt]
&&\quad  =2\int_{N_e(A)}\int_{-\infty}^{\infty}{
\Bbb E}\xi \bigl(\mathbf{0}_x+rn,\mathcal{P}_{\kappa(x)}^\mathrm{hom},H_x
\bigr) \kappa(x) \,\mathrm{d}r \nu _{d-1} \bigl(\,\mathrm{d}(x,n) \bigr).\nonumber
\end{eqnarray}

Now, we simplify the last integral and show that it coincides with
$\mu(\xi, \partial A)$, as given in \eqref{expec1}. First, recall
that there is a unique unit normal vector $n(x)$ at each $x \in
\partial^{1}A$ and define a measure $\mu_{d-1}$ on $N(A)$ by
\[
\mu_{d-1}( \cdot )={1\over2}\int_{\partial^{1}A}\mathbf{1} \bigl( \bigl(x,n(x) \bigr)\in \cdot \bigr) {\mathcal{H}}^{d-1}(
\mathrm{d}x).
\]
Since $A\in\mathbf{A}(d)$ it follows by Corollary~2.5 and Proposition~4.1 in
\cite{HLW} that
\[
\mu_{d-1}( \cdot )={1\over2}\int_{N(A)}\mathbf{1} \bigl((x,n)\in \cdot \bigr) H_0(x,n) {\mathcal{H}}^{d-1}
\bigl(\mathrm{d}(x,n) \bigr),
\]
where $H_0(x,n)$ is a certain function depending on the so-called
generalized principal curvatures of $A$; see equations (2.13) and
(2.24) in \cite{HLW}. Next, write
\begin{eqnarray*}
  \int_{N_e(A)} f(x,n)
\nu_{d-1} \bigl(\mathrm{d}(x,n) \bigr)&=&\int_{N(A)}f(x,n)
\mu _{d-1} \bigl(\mathrm{d}(x,n) \bigr)
\\
&&{}+\int_{T(N(A^*))}f(x,n) \mu_{d-1} \bigl(\mathrm{d} (x,n)
\bigr)
\\
&&{}-\int_{N(A)\cap T(N(A^*))}f(x,n) \mu_{d-1} \bigl(\mathrm{d}(x,n)
\bigr).   %
\end{eqnarray*}
According to the discussion before Theorem~5.2 in \cite{HLW}, given a
measurable function $f$ on $\mathbb{R}^d\times{\Bbb S}^{d-1}$, we can split
the integral over $N_e(A)$ in \eqref
{eq:AlmostFinishedExpectationAsymptotic} into three parts. The
projection map $\pi_1:N(A)\to{\Bbb R}^d, (x,n)\mapsto x$ has
Jacobian also given by $H_0(x,n)$ for ${\mathcal{H}}^{d-1}$-almost all
$(x,n)\in
N(A)$; see \cite{HLW}, Section~3. Combining these facts with the area
formula \cite{Federer}, Paragraph~3.2.3, applied to $\pi_1$ in each of
the three resulting integrals, which can be combined to a single
integral over $\partial^1 A$, we find that
\begin{eqnarray*}
  &&2\int_{N_e(A)}\int
_{-\infty}^{\infty}{\Bbb E}\xi \bigl(\mathbf{0}_x+rn,\mathcal{P}_{\kappa(x)}^\mathrm{hom},H_x
\bigr)\kappa(x) \,\mathrm{d}r \nu _{d-1} \bigl(\,\mathrm{d}(x,n) \bigr)
\\
&&\quad =\int_{\partial^1 A}\int_{-\infty}^\infty{
\Bbb E}\xi \bigl(\mathbf{0}_x+rn(x),\mathcal{P}_{\kappa(x)}^\mathrm{hom},H_x
\bigr) \kappa(x) H_0 \bigl(x,n(x) \bigr) H_0
\bigl(x,n(x) \bigr)^{-1} \,\mathrm{d}r \mathcal{H}^{d-1}(\mathrm {d}x)
\\
&&\quad =\int_{\partial^1 A}\int_{-\infty}^\infty{
\Bbb E}\xi \bigl(\mathbf{0}_x+rn(x),\mathcal{P}_{\kappa(x)}^\mathrm{hom},H_x
\bigr) \kappa(x) \,\mathrm{d}r \mathcal{H}^{d-1}(\mathrm{d}x),
\end{eqnarray*}
where we also have used the explicit representation of the measure $\mu
_{d-1}$ as well as the fact that ${\mathcal{H}}^{d-1}({\partial
A}^{++}\setminus
\partial^1A)=0$, which holds because $A\in\mathbf{A}(d)$. Since $\xi
$ is
invariant under rotations, we may replace $H_x$ by $\mathbb{R}^{d-1}$ and
$\mathbf{0}_x+rn(x)$ by $\mathbf{0}+rn$ to obtain \eqref{expec1} from \eqref
{eq:lim2}, as desired.
\end{pf*}

\begin{pf*}{Proof of Theorem~\ref{eq:Var}}
Applying the Mecke formula
for Poisson point processes, we get
\begin{eqnarray}
\label{Mecke} \lambda ^{-(d-1)/d} \operatorname{Var} \bigl[H^\xi({
\mathcal{P}}_\lambda, \partial A) \bigr] &=& \lambda ^{1/d} \int
_{\mathbb{R}^d} \mathbb{E}\xi_\lambda (x, {
\mathcal{P}}_\lambda, \partial A)^2\kappa (x) \,\mathrm{d}x\nonumber
\\[-8pt]\\[-8pt]
&&{}+ \lambda ^{1 + 1/d} \int_{\mathbb{R}^d} \int
_{\mathbb{R}^d} I_1 \kappa (x)\kappa (w) \,\mathrm{d}w
\,\mathrm{d}x,\nonumber
\end{eqnarray}
where
\[
I_1:= \mathbb{E}\xi_\lambda \bigl(x, {\mathcal{P}}_\lambda
\cup\{w\}, \partial A \bigr)\xi_\lambda \bigl(w, {\mathcal{P}}
_\lambda \cup\{x\}, \partial A \bigr) - \mathbb{E}\xi_\lambda (x,
{\mathcal{P}}_\lambda, \partial A) \mathbb{E}\xi_\lambda (w, {
\mathcal{P}}_\lambda, \partial A).
\]
The proof of Theorem~\ref{eq:WLLN} shows that the first integral in
\eqref{Mecke} converges to
\[
\int_{\partial^1 A}\int_{-\infty}^\infty{\Bbb
E}\xi \bigl(\mathbf{0}_x+rn(x),\mathcal{P}_{\kappa(x)}^\mathrm{hom},
\mathbb{R}^{d-1} \bigr)^2 \kappa(x) \,\mathrm{d}r
\mathcal{H}^{d-1}(\mathrm{d}x) = \mu \bigl(\xi^2, \partial A
\bigr).
\]
To complete the proof, we show that the second integral in \eqref
{Mecke} converges to the quadruple integral
in \eqref{eq: sigma}. We re-write the integral with respect to $x$
according to the generalized Steiner formula \eqref{eq:Steiner}, using
the notation already introduced in the proof of Theorem~\ref{eq:WLLN}.
Furthermore, for all $(x,n) \in N_e(A)$, let $H(x,n)$ denote the
hyperplane orthogonal to $n$ and containing $x$. Given $(x,n) \in
N_e(A)$, we re-write the integral with respect to $w$ as the iterated
integral over $H(x,n)$ and $\mathbb{R}$.
%$\int_{v \in H(x,n)} \int_{s \in\R}... dv ds$., and for $p \in H(x,n)$ and $r\in T(x,n)$, let $S(x,n,p,r):= \{s \in
%\R:  x + p + (r + s)n \in Q \}$.
This gives
\begin{eqnarray*}
&&\lambda ^{1 + 1/d} \int_{\mathbb{R}^d} \int
_{\mathbb{R}^d} I_1 \kappa (x)\kappa (w) \,\mathrm{d}x
\,\mathrm{d}w
\\
&&\quad = \lambda ^{1 + 1/d} \sum_{j = 1}^{d-1}
\omega_{d- j} \int_{(x,n)
\in
N_e(A)} \int_{r \in T(x,n)}
\int_{v \in H(x,n)} \int_{s \in\mathbb{R}} r^{d-1-j}
I_2
\\
&&\qquad {} \times\kappa (x + rn)\kappa \bigl((x + rn) + (v + sn) \bigr) \,\mathrm{d}s
\,\mathrm{d}v \,\mathrm{d}r \nu_j \bigl(\mathrm{d}(x,n) \bigr)
\end{eqnarray*}
with $I_2$ equal to
\begin{eqnarray*}
&&\mathbb{E}\xi_\lambda \bigl(x + rn, {\mathcal{P}}_\lambda \cup
\bigl\{(x + rn) + (v + sn) \bigr\}, \partial A \bigr)\xi_\lambda \bigl((x
+ rn) + (v + sn), {\mathcal{P}}_\lambda \cup\{ x + rn\}, \partial A
\bigr)
\\
&&\quad {} -\mathbb{E}\xi_\lambda (x + rn, {\mathcal{P}}_\lambda,
\partial A) \mathbb{E}\xi_\lambda \bigl((x + rn) + (v + sn), {
\mathcal{P}}_\lambda, \partial A \bigr).
\end{eqnarray*}
We change variables by putting $s = \lambda ^{-1/d}s'$, $r = \lambda
^{-1/d}r'$
and $v = \lambda ^{-1/d}v'$. This transforms the differential
$\lambda ^{1 + 1/d} \,\mathrm{d}s\,\mathrm{d}v\,\mathrm{d}r \nu_j(\mathrm
{d}(x,n))$ into
\[
\mathrm{d}s'\,\mathrm{d}v'\,\mathrm{d}r'
\nu_j \bigl(\mathrm{d}(x,n) \bigr),\qquad  j \in\{1,\ldots,d-1\}
\]
and $I_2$ into $I_3$ given by
\begin{eqnarray*}
I_3&:=& \mathbb{E}\xi_\lambda \bigl(x + \lambda
^{-1/d}r'n, {\mathcal {P}}_\lambda \cup \bigl\{
\bigl(x + \lambda ^{-1/d}r'n \bigr) + \bigl(\lambda
^{-1/d}v' + \lambda ^{-1/d}s'n \bigr)
\bigr\}, \partial A \bigr)
\\
&&{} \times\xi_\lambda \bigl( \bigl(x + \lambda ^{-1/d}r'n
\bigr) + \bigl(\lambda ^{-1/d}v' + \lambda
^{-1/d}s'n \bigr), {\mathcal{P}}_\lambda \cup \bigl
\{ x + \lambda ^{-1/d}r'n \bigr\}, \partial A \bigr)
\\
&&{}-\mathbb{E}\xi_\lambda \bigl(x + \lambda ^{-1/d}r'n,
{\mathcal {P}}_\lambda, \partial A \bigr) \mathbb{E}\xi_\lambda
\bigl( \bigl(x + \lambda ^{-1/d}r'n \bigr) + \bigl(\lambda
^{-1/d}v' + \lambda ^{-1/d}s'n
\bigr), {\mathcal{P}}_\lambda, \partial A \bigr).
\end{eqnarray*}
To simplify the notation, we shall write $s$, $r$ and $v$ for $s'$,
$r'$ and $v'$, respectively. Then
\begin{eqnarray}
\label{sum}
  &&\lambda ^{1 + 1/d} \int
_{\mathbb{R}^d} \int_{\mathbb{R}^d} I_1 \kappa
(x)\kappa (w) \,\mathrm{d} x \,\mathrm{d}w\nonumber
\\
&&\quad  = \sum_{j = 1}^{d-1} \lambda ^{-(d - 1 - j)/d}
\omega_{d- j} \int_{N_e(A)} \int_{\lambda ^{1/d} T(x,n)}
\int_{H(x,n)} \int_{\mathbb{R}} r^{d-1-j}
I_3
\\
&&\qquad {} \times\kappa \bigl(x + \lambda ^{-1/d}rn \bigr)\kappa \bigl( \bigl(x +
\lambda ^{-1/d}rn \bigr) + \bigl(\lambda ^{-1/d}v + \lambda
^{-1/d}sn \bigr) \bigr) \,\mathrm{d}s\,\mathrm{d}v\,\mathrm{d}r\nu
_j \bigl(\mathrm{d}(x,n) \bigr). \nonumber  %
\end{eqnarray}
By the second part of Lemma~\ref{eq:L3}, the factor $\llvert  I_3\rrvert   $ in \eqref
{sum} is dominated
uniformly in $\lambda $ by an integrable function of $(x,n) \in
N_e(A)$, $s
\in\mathbb{R}$, $v\in H(x,n)$ and $r \in\mathbb{R}$. More precisely,
%there is a constant $C \in(0, \infty)$ such that for all $\la\in[1,
%\infty)$, %$p \in\HH(x,n)$, $u \in\R$, and
%$s \in\R$ we have
\[
\llvert I_3\rrvert \leq c G^{\xi,p} \bigl(\llvert r\rrvert
\bigr)^{1/p} \exp \bigl( -c^{-1} \sqrt{\llVert v\rrVert
^2 + s^2} \bigr),
\]
where the constant $c$ is independent of all arguments. Thus for each
$j \in\{1,\ldots,d-1\}$, we have
\begin{eqnarray*}
&& \left\vert \int_{N_e(A)} \int_{\lambda ^{1/d} T(x,n)} \int
_{H(x,n)} \int_{\mathbb{R}} r^{d-1-j}\right.
I_3
\\
&&\qquad {} \times\left.\kappa \bigl(x + \lambda ^{-1/d}rn \bigr)\kappa \bigl( \bigl(x +
\lambda ^{-1/d}rn \bigr) + \bigl(\lambda ^{-1/d}v + \lambda
^{-1/d}sn \bigr) \bigr) \,\mathrm{d}s\,\mathrm{d}v\,\mathrm{d}r\nu
_j \bigl(\mathrm{d}(x,n) \bigr) \vphantom{ \int_{N_e(A)}}\right\vert
\\
&&\quad \leq c\llVert \kappa \rrVert _\infty^2\int
_{N_e(A)}\int_{-\infty}^{\infty
}\int
_{\mathbb{R}^{d-1}}\int_{\mathbb{R}} r^{d-1-j}
G^{\xi,p} \bigl(\llvert r\rrvert \bigr)^{1/p}
\\
&&\qquad {}\times  \exp \bigl( -
c^{-1} \sqrt{\llVert v\rrVert ^2 + s^2} \bigr)
\,\mathrm{d}s\,\mathrm{d}v\,\mathrm{d}r\nu_j \bigl(\mathrm {d}(x,n) \bigr)
\\
&&\quad \leq c\llVert \kappa \rrVert _\infty^2\llvert
\nu_j\rrvert \bigl(N_e(A) \bigr) \int
_{-\infty
}^{\infty
}\int_{\mathbb{R}^{d-1}}\int
_{\mathbb{R}} r^{d-1-j} G^{\xi,p} \bigl(\llvert r
\rrvert \bigr)^{1/p} \exp \bigl( - c^{-1} \sqrt{\llVert v
\rrVert ^2 + s^2} \bigr) \,\mathrm{d}s\,\mathrm{d}v\,\mathrm{d}r.
\end{eqnarray*}
Notice that $\llvert  \nu_j\rrvert   (N_e(A))$ and the triple integral are finite by
the assumption that $A\in\mathbf{A}(d)$ and the moment condition
\eqref
{eq:mom}, respectively.
As in the proof of Theorem~\ref{eq:WLLN}, we have
$\lim_{\lambda \to\infty}\lambda ^{1/d} T(x,n) = (-\infty, \infty
)$. %$ \mbox{and} \liml\la^{1/d} S(x,n,\la^{-1/d}p,\la^{-1/d}s
%) = (0, \infty).$$
Taking the limit, as $\lambda \to\infty$, in \eqref{sum} and applying
the dominated convergence theorem,
we see that only the term $j = d - 1$ remains. By Fubini's theorem and
Lemma~\ref{L2}, this gives
\begin{eqnarray*}
&&\lim_{\lambda\to\infty} \lambda ^{1 + 1/d} \int_{\mathbb{R}^d}
\int_{\mathbb{R}^d} I_1 \kappa (x)\kappa (w) \,\mathrm{d}x
\,\mathrm{d}w
\\
&&\quad  = 2\int_{N_e(A)}\int_{\mathbb{R}^{d-1}}\int
_{-\infty}^\infty \int_{-\infty}^\infty
c^\xi \bigl(\mathbf{0}_x+rn,v+sn;\mathcal{P}^{\mathrm
{hom}}_{\kappa (x)},
\mathbb{R}^{d-1} \bigr) \kappa(x)^2 \,\mathrm{d}s\,\mathrm{d}r
\,\mathrm{d}v\nu_{d-1} \bigl(\mathrm {d}(x,n) \bigr).
\end{eqnarray*}
We can now use the same arguments as in the proof of Theorem~\ref
{eq:WLLN} to show that the integral reduces to the
quadruple integral in \eqref{eq: sigma}. This yields \eqref{eq:Var1},
as desired.
\end{pf*}

%s4
\section{Proofs of Theorems \texorpdfstring{\protect\ref{thm:VoronoiV+S}}{1.1}--\texorpdfstring{\protect\ref{thm:IteratedVoronoi}}{1.5}}\label{sec:ProofsApplications}

We shall deduce Theorems \ref{thm:VoronoiV+S}--\ref
{thm:IteratedVoronoi} from the general Theorems \ref{eq:WLLN} and
\ref{eq:Var}. In each case, it suffices to express the relevant
statistic as a sum of score functions and to show that the score function
satisfies the conditions of the general theorems.
We anticipate that the expectation formula \eqref{expec-formula} could
be evaluated explicitly for some of the score functions described below.
The proof of the positivity of the constants appearing in the variance
expressions is postponed to Section~\ref{sec:VLB}.

\begin{pf*}{Proof of Theorem~\ref{thm:VoronoiV+S}}
We first prove the asserted results for the volume functional
$V_\lambda
(A)$, with $A\in\mathbf{A}(d)$.
For locally finite $\mathcal{X}\subset\mathbb{R}^d$, $x \in\mathcal
{X}$, define the score
\begin{equation}
\label{defp} \xi^{(1)}(x, \mathcal{X}, \partial A):= \cases{ \operatorname{Vol} \bigl(v(x) \cap A^c
\bigr), &\quad  if  $v(x) \cap \partial A \neq\varnothing$, $x \in A$,
\cr
- \operatorname{Vol} \bigl(v(x) \cap A \bigr), &\quad  if  $v(x) \cap\partial A
\neq\varnothing$, $x \in A^c$,
\cr
0, &\quad if  $v(x) \cap\partial A = \varnothing $, }
\end{equation}
where $v(x):=v(x, \mathcal{X})$ is the Voronoi cell of $x$ based on
the point
configuration $\mathcal{X}$. In view of the limits appearing in our main
results, we also need to define scores on hyperplanes, that is, on~$\mathbb{R}
^{d-1}$. We thus put
\begin{equation}
\label{defpinf} \xi^{(1)} \bigl(x, \mathcal{X}, \mathbb{R}^{d-1}
\bigr):= \cases{
 \operatorname{Vol} \bigl(v(x)
\cap \mathbb{R}_+^{d-1} \bigr), &\quad  if  $v(x) \cap
\mathbb{R}^{d-1} \neq\varnothing$, $x \in\mathbb{R}_{-}^{d-1}$,
\cr
- \operatorname{Vol} \bigl(v(x) \cap\mathbb{R}_{-}^{d-1}
\bigr), &\quad  if  $v(x) \cap\mathbb{R}^{d-1} \neq\varnothing$, $x \in
\mathbb{R}_+^{d-1}$,
\cr
0, & \quad if  $v(x) \cap\mathbb{R}^{d-1} = \varnothing $, }
\end{equation}
where we recall $\mathbb{R}_+^{d-1}:= \mathbb{R}^{d-1} \times[0,
\infty)$ and
$\mathbb{R}_{-}^{d-1}:= \mathbb{R}^{d-1} \times(-\infty, 0]$.
These definitions ensure that
\[
V_\lambda (A) - \operatorname{Vol}(A) = \sum
_{x \in{\mathcal
{P}}_\lambda } \xi^{(1)}(x, {\mathcal{P}}_\lambda,
\partial A) = \lambda ^{-1} \sum_{x \in{\mathcal{P}}_\lambda } \xi
^{(1)}_\lambda (x, {\mathcal{P}} _\lambda, \partial A),
\]
where we use that $\xi^{(1)}$ is homogenous of order $d$.
We wish to deduce the volume asymptotics for $ V_\lambda (A)$ by applying
the limits \eqref{eq: newexpect} and \eqref{eq: newvar} with $\gamma
= d$ and
with $\xi$ set to $\xi^{(1)}$.
It suffices to show that the score $\xi^{(1)}$ is homogeneously
stabilizing \eqref{eq:hom}, exponentially stabilizing as at \eqref
{eq:expo}, satisfies the moment condition \eqref{eq:mom} for $p=1$ and
some $p \in(2, \infty)$, and is well approximated by ${\mathcal
{P}}_\lambda $ input
on half-spaces as at \eqref{eq:lin}.
The first three conditions have been established several times in the
literature; see the proof of Theorem~2.2 of \cite{Yu}.
%We wish to deduce the volume asymptotics for $ V_\la(A)$ by applying
%the limits \eqref{eq: newexpect} and \eqref{eq: newvar} with $\gamma=
%d$.

To show that $\xi^{(1)}$ is well approximated by ${\mathcal
{P}}_\lambda $ input on
half-spaces as at \eqref{eq:lin}, it suffices to slightly modify the proof
of the analogous result in Theorem~2.2 of \cite{Yu}. For the sake of
completeness, we provide the details as follows.

By definition of $\mathbf{A}(d)$, almost all points of $\partial A$
belong to
$\partial^{1} A$
and it so suffices to show \eqref{eq:lin} for a fixed $y \in\partial
^{1} A$.
Translating $y$ to the origin, letting ${\mathcal{P}}_\lambda $
denote a Poisson
point process on~$\mathbb{R}^d$, letting $\partial A$ denote
$\partial A - y$, and using rotation invariance of $\xi^{(1)}$, it is
enough to show for all $w \in\mathbb{R}^d$ that
\[
\lim_{\lambda \to\infty}\mathbb{E} \bigl\llvert \xi^{(1)} \bigl(w,
\lambda ^{1/d}{\mathcal{P}}_\lambda, \lambda ^{1/d}
\partial A \bigr) - \xi^{(1)} \bigl(w, \lambda ^{1/d}{
\mathcal{P}}_\lambda, \mathbb{R}^{d-1} \bigr) \bigr\rrvert = 0,
\]
where $\mathbb{R}^{d-1}$ is the unique hyperplane tangent to $\partial
A$ at
the origin.
Without loss of generality, we
assume, locally around the origin, that $\partial A \subset
\mathbb{R}_{-}^{d-1}$. Fix $\varepsilon> 0$ and $w \in\mathbb
{R}^d$. We note that
there is a constant $L \in(0, \infty)$ such that
\[
\sup_{\lambda >0} \bigl(\mathbb{E} \bigl[\xi^{(1)} \bigl(w,
\lambda ^{1/d} {\mathcal{P}}_\lambda, \lambda ^{1/d}
\partial A \bigr)^2 \bigr] \bigr)^{1/2} \leq L
\]
and
\[
\sup_{\lambda >0} \bigl(\mathbb{E} \bigl[\xi^{(1)} \bigl(w,
\lambda ^{1/d} {\mathcal{P}}_\lambda, \mathbb{R} ^{d-1}
\bigr)^2 \bigr] \bigr)^{1/2} \leq L.
\]
Let $\widetilde{v}(w,\lambda ^{1/d}{\mathcal{P}}_\lambda )$ be the
union of $v(w,\lambda
^{1/d}{\mathcal{P}}_\lambda )$ and all the Voronoi cells adjacent to
$v(w,\lambda ^{1/d}{\mathcal{P}}_\lambda )$ in the Voronoi mosaic of
${\mathcal{P}}_\lambda $. For all
$r \in(1, \infty)$,
consider the event
\begin{equation}
\label{eq:defE}E_1(\lambda,w,r):= \bigl\{ \operatorname{diam} \bigl(
\widetilde{v} \bigl(w, \lambda ^{1/d}{\mathcal{P}}_\lambda \bigr)
\bigr) \leq r \bigr\},
\end{equation}
where $\operatorname{diam}( \cdot )$
stands for the diameter of the argument set. Lemma~2.2 of \cite{MY}
shows there is $r_0:=r_0(\varepsilon, L)$ such that for $r \in[r_0,
\infty)$
and $\lambda $ large we have
$\mathbb{P}(E_1(\lambda,w,r)^c) \leq( \varepsilon/2L)^2$.
%Since $\xi^{(1)}(w, \la^{1/d}\P_\la, \la^{1/d}\partial A)$ and $
%\xi^{(1)}(w,
%\la^{1/d}\P_\la, \R^{d-1})$ have finite second moments, uniformly in
%$w \in\R^d$
%and $\la\in(0, \infty)$,
It follows by the Cauchy--Schwarz inequality that
\[
\lim_{\lambda \to\infty}\mathbb{E} \bigl\llvert \bigl(\xi^{(1)}
\bigl(w, \lambda ^{1/d}{\mathcal{P}}_\lambda, \lambda
^{1/d}\partial A \bigr) - \xi^{(1)} \bigl(w, \lambda
^{1/d}{\mathcal{P}}_\lambda, \mathbb {R}^{d-1} \bigr)
\bigr) \mathbf{1} \bigl(E_1(\lambda,w,r_0)^c
\bigr) \bigr\rrvert \leq\varepsilon.
\]
By the triangle inequality and the arbitrariness of $\varepsilon$, it
is therefore enough to show that
\begin{equation}
\label{eq:diff} \lim_{\lambda \to\infty} \mathbb{E} \bigl\llvert \bigl(\xi
^{(1)} \bigl(w, \lambda ^{1/d}{\mathcal{P}}_\lambda,
\lambda ^{1/d}\partial A \bigr) - \xi^{(1)} \bigl(w, \lambda
^{1/d}{\mathcal {P}}_\lambda, \mathbb{R}^{d-1} \bigr)
\bigr) \mathbf{1} \bigl(E_1(\lambda,w,r_0) \bigr) \bigr
\rrvert \leq\varepsilon.
\end{equation}
By the way that $y$ was chosen, $\mathbf{0}$ is a point in $\partial^{1}A$
and thus has a unique normal vector. We first assume $w \in\mathbb
{R}_{-}^{d-1}$;
the arguments with $w \in\mathbb{R}_{+}^{d-1}$ are nearly identical.
Moreover, we may assume $w \in
\lambda ^{1/d}A$ for $\lambda $ large.
Consider the (possibly degenerate) solid
\begin{equation}
\label{eq:Delta} \Delta_\lambda (w):= \Delta_\lambda
(w,r_0):= \bigl(\mathbb {R}^{d-1}_{-} \setminus
\lambda ^{1/d}A \bigr) \cap B_{ r_0 }(w).
\end{equation}
Recalling that $\partial A$ is $(\mathcal{H}^{d-1}, d - 1)$
rectifiable, it
follows that almost all of $\partial A$ is contained in a union of
$C^1$ sub-manifolds of $\mathbb{R}^d$ \cite{Federer}, Theorem~3.2.29.
Since $\mathbf{0}$ is a point of $\partial^{1}A$, it follows that %
%and since
%the boundary of $A$ is rectifiable,
the maximal `height' $h_\lambda (w, r_0)$ of the
solid $\Delta_\lambda (w,r_0)$ with respect to the hyperplane
$\mathbb{R}^{d-1}$ satisfies $\lim_{\lambda \to\infty}
h_\lambda (w, r_0)=
0$ for fixed $w$ and $r_0$ (see also the linear approximation
properties of rectifiable sets summarized in Chapter~15 of \cite
{Mattila}). Hence,
\[
\operatorname{Vol} \bigl(\Delta_\lambda (w,r_0) \bigr) = O
\bigl( h_\lambda (w, r_0) \cdot r_0^{d-1}
\bigr)
\]
%
%$$\Vol(\Delta_\la(w)) = O((\left\Vert w\right\Vert   + 2\b\log\la) \la^{-1/d} (2\b\log
%\la)^{d-1}) = O((\log\la)^d \la^{-1/d}).
%$$
and so for large $\lambda $ we have $\operatorname{Vol}(\Delta
_\lambda (w,r_0)) \leq
\varepsilon$.
On the event $E_1(\lambda,w,r_0)$, the difference of the volumes $v(w,
\lambda ^{1/d}{\mathcal{P}}_\lambda ) \cap\lambda ^{1/d}A^c$ and
$v(w, \lambda ^{1/d}{\mathcal{P}}_\lambda )
\cap
\mathbb{R}_+^{d-1}$
is at most
$\operatorname{Vol}(\Delta_\lambda (w,r_0))$. Thus, for large
$\lambda $ we get
\begin{eqnarray*}
&&\mathbb{E} \bigl\llvert \bigl(\xi^{(1)} \bigl(w, \lambda
^{1/d}{ \mathcal{P}}_\lambda, \lambda ^{1/d}\partial A
\bigr) - \xi^{(1)} \bigl(w, \lambda ^{1/d}{
\mathcal{P}}_\lambda, \mathbb{R}^{d-1} \bigr) \bigr) \mathbf{1}
\bigl(E_1(\lambda,w,r_0) \bigr) \bigr\rrvert
\\
&&\quad \leq\operatorname{Vol} \bigl(\Delta_\lambda (w,r_0) \bigr)
\leq\varepsilon,
\end{eqnarray*}
which gives \eqref{eq:lin} as desired. %[ If this is correct then I
%need to change the notation and proofs in the sequel. JY]

We now prove the asserted results for
the surface area functional $S_\lambda (A)$. As in \cite{Yu}, given
$\mathcal{X}$
locally finite and an admissible set $A \subset\mathbb{R}^d$, define
for $x
\in\mathcal{X}\cap A$
the area score $\xi^{(2)}(x, \mathcal{X}, \partial A)$ to be the
${\mathcal{H}}
^{d-1}$-measure of the $(d-1)$-dimensional faces of $v(x)$ belonging to the
boundary of $\bigcup_{x \in\mathcal{X}\cap A} v(x)$ (if there are no such
faces or
if $x \notin\mathcal{X}\cap A$, then put $\xi^{(2)}(x, \mathcal{X},
\partial A)$ to
be zero).
Similarly, for $x \in\mathcal{X}\cap\mathbb{R}^{d-1}_{-}$, put $\xi
^{(2)}(x, \mathcal{X},
\mathbb{R}^{d-1})$ to be
the ${\mathcal{H}}^{d-1}$-measure of the $(d-1)$-dimensional faces of $v(x)$
belonging to the\vspace*{-1pt}
boundary of $\bigcup_{x\in\mathcal{X}\cap\mathbb{R}^{d-1}_{-}} v(x)$,
otherwise $\xi^{(2)}(x, \mathcal{X}, \mathbb{R}^{d-1})$ is zero. We
note that $\xi
^{(2)}$ is homogenous of order $d - 1$ and that
%\CT[Slightly modifications here..]
\[
S_\lambda (A) = \sum_{x \in{\mathcal{P}}_\lambda }
\xi^{(2)}(x, {\mathcal{P}}_\lambda, {\partial A}).
\]
We wish to deduce the first- and second-order limit behavior of
$S_\lambda
(A)$ by applying the limits \eqref{eq: newexpect} and \eqref{eq:
newvar} with $\gamma= d-1$ and with $\xi$ set to $\xi^{(2)}$.

It is easy to see and well known that the score $\xi^{(2)}$ is
homogeneously stabilizing \eqref{eq:hom}, exponentially stabilizing
\eqref{eq:expo}, and satisfies the moment condition \eqref{eq:mom}
for all $p \geq1$; see, for example, the proof of Theorem~2.4 of \cite{Yu}.
To see that $\xi^{(2)}$ is well approximated by ${\mathcal
{P}}_\lambda $ input on
half-spaces \eqref{eq:lin}, it suffices to follow the proof of Theorem~2.4 of \cite{Yu}. For sake of completeness, we include the details as follows.

Fix $\varepsilon> 0$ and $w \in\mathbb{R}^d$. By the moment bounds
on $\xi
^{(2)}$ and the
Cauchy--Schwarz inequality, it is enough to show the following
counterpart to
\eqref{eq:diff}, namely to show that
\begin{equation}
\label{eq:adiff} \lim_{\lambda \to\infty} \mathbb{E} \bigl\llvert \bigl(
\xi^{(2)} \bigl(w, \lambda ^{1/d}{\mathcal{P}}_\lambda,
\lambda ^{1/d}\partial A \bigr) - \xi^{(2)} \bigl(w, \lambda
^{1/d}{\mathcal {P}}_\lambda, \mathbb{R}^{d-1} \bigr)
\bigr) \mathbf{1} \bigl(E_1(\lambda,w,r_0) \bigr) \bigr
\rrvert \leq C \varepsilon^{1/2},
\end{equation}
where
$E_1(\lambda,w,r_0)$ is as at \eqref{eq:defE}, and where, as above,
the origin is a point of $\partial^1A-y$.
Define
\[
E_0(\lambda, w,r_0):= \bigl\{ \lambda ^{1/d}
{\mathcal{P}}_\lambda \cap \Delta_\lambda (w,r_0) =
\varnothing \bigr\},
\]
where $\Delta_\lambda (w,r_0)$ is as at \eqref{eq:Delta}.
The intensity measure of $\lambda ^{1/d} {\mathcal{P}}_\lambda $ is
upper bounded by
$\llVert  \kappa \rrVert   _{\infty}$, yielding
for large $\lambda $ that
\begin{equation}
\label{eq:EE} \mathbb{P} \bigl[E_0(\lambda, w,r_0)^c
\bigr] \leq1 - \exp \bigl(-\llVert \kappa \rrVert _{\infty}
\operatorname{Vol} \bigl(\Delta_\lambda (w,r_0) \bigr) \bigr)
\leq\llVert \kappa \rrVert _{\infty} \varepsilon,
\end{equation}
%
%$$ \leq1 - \exp(- c_6 (\log\la)^d \la^{-1/d} ) = O((\log\la)^d
%\la^{-1/d}).$$
where we used that $\operatorname{Vol}(\Delta_\lambda (w,r_0)) \leq
\varepsilon$.

The two score functions $\xi^{(2)}(w,
\lambda ^{1/d}{\mathcal{P}}_\lambda, \lambda ^{1/d}\partial A)$
and $\xi^{(2)}(w, \lambda
^{1/d}{\mathcal{P}}_\lambda,
\mathbb{R}^{d-1})$ coincide on the event $E_1(\lambda, w,r_0) \cap
E_0(\lambda,
w,r_0)$. Indeed, on this event it follows that
$f$ is a face of a boundary cell of $\lambda ^{1/d}A_\lambda $ iff
$f$ is a face
of a boundary cell of the Poisson--Voronoi mosaic of
$\mathbb{R}_{-}^{d-1}$. (If $f$ is a face of the boundary cell $v(w,
\lambda ^{1/d} {\mathcal{P}}_\lambda ), w \in\lambda ^{1/d}A$,
then $f$ is also a face of
$v(z, \lambda ^{1/d} {\mathcal{P}}_\lambda )$ for some $z \in
\lambda ^{1/d}A^c$. If
$\lambda ^{1/d} {\mathcal{P}}_\lambda \cap\Delta_\lambda (w,r_0)
= \varnothing$, then $z$
must belong to $\mathbb{R}_{+}^{d-1}$,
showing that $f$ is face of a boundary cell of the Poisson--Voronoi
mosaic of $\mathbb{R}_{-}^{d-1}$. The reverse implication is shown
similarly.)

On the other hand, since
%\CT[Is the square at the right place here?]
\[
\mathbb{E} \bigl[ \bigl(\xi^{(2)} \bigl(w, \lambda ^{1/d}{
\mathcal{P}}_\lambda, \lambda ^{1/d}\partial A \bigr) -
\xi^{(2)} \bigl(w, \lambda ^{1/d}{\mathcal{P}}_\lambda,
\mathbb{R}^{d-1} \bigr) \bigr)^2 \mathbf{1}
\bigl(E_1( \lambda, w,r_0) \bigr) \bigr] = O(1),
\]
and since by \eqref{eq:EE} we have $\mathbb{P}[E_{0}(\lambda, w,r_0)^c]
\leq\llVert  \kappa \rrVert   _{\infty} \varepsilon$, it follows by the Cauchy--Schwarz
inequality that, as $\lambda \to\infty$,
\begin{eqnarray}
\label{eq:diff-alpha-2}
  && \mathbb{E} \bigl\llvert \bigl(
\xi^{(2)} \bigl(w, \lambda ^{1/d}{\mathcal{P}}_\lambda,
\lambda ^{1/d}\partial A \bigr) - \xi^{(2)} \bigl(w, \lambda
^{1/d}{\mathcal {P}}_\lambda, \mathbb{R}^{d-1} \bigr)
\bigr) \mathbf{1} \bigl(E_1(\lambda,w,r_0) \bigr) \bigr
\rrvert\nonumber
\\
&&\quad  = \mathbb{E}\bigl\vert \bigl(\xi^{(2)} \bigl(w, \lambda ^{1/d}{
\mathcal{P}}_\lambda, \lambda ^{1/d}\partial A \bigr) -
\xi^{(2)} \bigl(w, \lambda ^{1/d}{\mathcal{P}}_\lambda,
\mathbb{R}^{d-1} \bigr) \bigr)\nonumber
\\[-8pt]\\[-8pt]
&&\qquad {} \times\mathbf{1} \bigl(E_1(\lambda, w,r_0) \bigr) \mathbf{1} \bigl(E_{0}(\lambda, w,r_0)^c \bigr)
\bigr\vert\nonumber
\\
&&\quad  \leq C \bigl(\llVert \kappa \rrVert _{\infty} \varepsilon
\bigr)^{1/2}. \nonumber  %
\end{eqnarray}
Therefore, \eqref{eq:adiff} holds and so $\xi^{(2)}$ is
well approximated by ${\mathcal{P}}_\lambda $ input on half-spaces as
at \eqref
{eq:lin}, as desired.
\end{pf*}

\begin{pf*}{Proof of Theorem~\ref{thm:VoronoiSkeleton}}
Let us first recall that the Poisson--Voronoi mosaic is a normal
mosaic; see \cite{SW}. This means that with probability one each $\ell
$-dimensional face in $\operatorname{skel}_\ell(\operatorname{PV}_\lambda (A))$ arises as
the intersection of exactly $d-\ell+1$ Voronoi cells.

Now, given $\mathcal{X}$ locally finite, $x\in\mathcal{X}$, and an admissible
$A\subset\mathbb{R}^d$, define $\xi^{(3,\ell)}(x,\mathcal
{X},\partial A)$ as
\[
\xi^{(3,\ell)}(x,\mathcal{X},\partial A):={1\over d-\ell}\mathop{\sum_{f\in{\mathcal{F}}_\ell(v(x))}}_{ f\subset\partial({\operatorname{PV}}_\lambda(A))}\mathcal{H}^{\ell}(f)
\]
%\[
%\xi^{(3,\ell)}(x,\mathcal{X},\partial A):={1\over d-\ell+1}\mathop{\sum
%_{f\in\mathcal{F}_\ell(v(x))}}_{ f\subset\partial(\operatorname{PV}_\lambda(A))}{\mathcal{H}} ^{(\ell)}(f)
%\]
%
and zero otherwise. Then
\[
H^{(\ell)}_\lambda (A) =\mathop{\sum_{x\in{\mathcal{P}}_\lambda }}_{
x\in A}
\xi^{(3,\ell
)}(x,{\mathcal{P}}_\lambda,\partial A) = \lambda
^{- \ell/d} \sum_{x\in{\mathcal{P}}_\lambda }\xi ^{(3,\ell)}_\lambda
(x,{\mathcal{P}}_\lambda,\partial A),
\]
where we used that $\xi^{(3,\ell)}$ is homogeneous of order $\ell$.
We wish to deduce the first- and second-order limit behaviour of
$H^{(\ell)}_\lambda (A)$ by applying the limits \eqref{eq:
newexpect} and
\eqref{eq: newvar} with $\gamma= \ell$ and with $\xi$ set to $\xi
^{(3,\ell)}$.

The proof that $\xi^{(3,\ell)}$ is homogeneously stabilizing \eqref
{eq:hom}, exponentially stabilizing \eqref{eq:expo}, and satisfies the
moment condition \eqref{eq:mom} for all $p \geq1$ follows nearly
verbatim the proof that $\xi^{(2)}$ has these properties. Indeed the
radius of stabilization for $\xi^{(3,\ell)}$ coincides with that of
$\xi^{(2)}$.

To see that $\xi^{(3,\ell)}$ is well approximated by ${\mathcal
{P}}_\lambda $ input
on half-spaces as at \eqref{eq:lin}, we may follow the proof that
$\xi^{(2)}$ is well approximated by ${\mathcal{P}}_\lambda $ input
on half-spaces.
Notice that on the event $E_1(\lambda, w,r_0) \cap E_0(\lambda,
w,r_0)$, the
scores $\xi^{(3,\ell)}(w,
\lambda ^{1/d}{\mathcal{P}}_\lambda, \lambda ^{1/d}\partial A)$
and $\xi^{(3,\ell)}(w, \lambda
^{1/d}{\mathcal{P}}_\lambda,
\mathbb{R}^{d-1})$ coincide. As in \eqref{eq:diff-alpha-2}, we obtain
\begin{eqnarray*}
& &\mathbb{E} \bigl\llvert \bigl(\xi^{(3,\ell)} \bigl(w, \lambda
^{1/d}{\mathcal {P}}_\lambda, \lambda ^{1/d}\partial A
\bigr) - \xi^{(3,\ell)} \bigl(w, \lambda ^{1/d}{
\mathcal{P}}_\lambda, \mathbb{R}^{d-1} \bigr) \bigr) \mathbf{1}
\bigl(E_1(\lambda, w,r_0) \bigr) \bigr\rrvert
\\
&&\quad  = \mathbb{E}\bigl\vert \bigl(\xi^{(3,\ell)} \bigl(w, \lambda ^{1/d}{
\mathcal{P}}_\lambda, \lambda ^{1/d}\partial A \bigr) -
\xi^{(3,\ell)} \bigl(w, \lambda ^{1/d}{\mathcal{P}}_\lambda,
\mathbb{R}^{d-1} \bigr) \bigr)
\\
&&\qquad {} \times\mathbf{1} \bigl(E_1(\lambda, w,r_0) \bigr) \mathbf{1} \bigl(E_{0}(\lambda, w,r_0)^c \bigr)
\bigr\vert \leq C \bigl(\llVert \kappa \rrVert _\infty\varepsilon
\bigr)^{1/2}.
\end{eqnarray*}
This gives that $\xi^{(3,\ell)}$ satisfies \eqref{eq:lin} as
desired.
\end{pf*}

\begin{pf*}{Proof of Theorem~\ref{thm:VoronoiCombinatorics}}
Given $\mathcal{X}$ locally finite, $x \in\mathcal{X}$, and $A\in
\mathbf{A}(d)$, let us
define the score
$\xi^{(4,\ell)}(x, \mathcal{X}, \partial A)$ to be the number of
$\ell
$-dimensional faces of $v(x):=v(x, \mathcal{X})$ belonging to
$\partial(\operatorname{PV}_\lambda(A))$. Define $\xi^{(4,\ell)}(x, \mathcal{X}, \mathbb
{R}^{d-1})$ similarly.
Then
\[
f_\lambda ^{\ell}(A) = \sum_{x \in{\mathcal{P}}_\lambda }
\xi ^{(4,\ell)}_\lambda (x, {\mathcal{P}} _\lambda, \partial
A).
\]
We shall show that $\xi^{(4,\ell)}$ satisfies the hypotheses of
Theorems \ref{eq:WLLN} and \ref{eq:Var}, and thus
deduce Theorem~\ref{thm:VoronoiCombinatorics} from \eqref{eq:
newexpect} and \eqref{eq: newvar} with $\xi$ set to
$\xi^{(4,\ell)}$ and $\gamma$ set to zero (notice that $\xi
^{(4,\ell)}$ is homogeneous of order $0$).
For brevity, write $\xi^{(4)}$ for $\xi^{(4,\ell)}$ for fixed $\ell
\in\{0,\ldots,d-1\}$. Now, $\xi^{(4)}$ is homogeneously and
exponentially stabilizing since its radius of stabilization coincides
with that for the
volume score $\xi^{(1)}$ defined in the proof of Theorem~\ref{thm:VoronoiV+S}.
The number $N^{(\ell)}(x, {\mathcal{P}}_\lambda )$ of $\ell
$-dimensional faces of a
Poisson--Voronoi cell $v(x)$ has moments of all orders and, therefore,
the moment condition \eqref{eq:mom} holds because
\begin{eqnarray*}
&& \bigl\llvert \xi^{(4)}_\lambda \bigl( x+\lambda
^{-1/d}rn, {\mathcal{P}}_\lambda \cup\{z\}, \partial A \bigr)
\bigr\rrvert
\\
&&\quad \leq N^{(\ell)} \bigl(x+\lambda ^{-1/d}rn, {\mathcal{P}}_\lambda
\cup\{ z\} \bigr) \mathbf{1} \bigl( v \bigl(\lambda ^{1/d}x+rn, \lambda
^{1/d} {\mathcal{P}}_\lambda \bigr) \cap \partial A \neq
\varnothing \bigr)
\end{eqnarray*}
for $(x,n)\in N_e(A)$.
The expectation of the last factor decays uniformly fast in $r$, giving
that $\xi^{(4)}$ satisfies the moment condition \eqref{eq:mom}
for all $p \geq1$.

The arguments in the proof of
Theorem~\ref{thm:VoronoiV+S} showing that the surface area score $\xi
^{(2)}$ is well approximated by ${\mathcal{P}}_\lambda $ input on
half-spaces extend
to show that
$\xi^{(4)}$ is likewise well approximated by ${\mathcal{P}}_\lambda
$ input on
half-spaces. %Theorem~2.4 of \cite{Yu}.
The guiding idea is that with high probability, we have that $f$ is a
face of a Voronoi cell $v(w)$ belonging to the Poisson--Voronoi
approximation of $\lambda ^{1/d}(A - y)$ if and only if it belongs to the
Poisson--Voronoi approximation of $\mathbb{R}^{d-1}_+$. Indeed, this happens
on the high probability event that the region `between' the boundary of
the Poisson--Voronoi approximation of $A$ and $\mathbb{R}^{d-1}$ in the
neighbourhood of $w$, must be devoid of points; see the proof of
Theorem~\ref{thm:VoronoiV+S}. Thus,
$\xi^{(4,\ell)}$ satisfies all the hypotheses of Theorems \ref
{eq:WLLN} and \ref{eq:Var} and this completes the proof of Theorem~\ref{thm:VoronoiCombinatorics}.
\end{pf*}

\begin{pf*}{Proof of Theorem~\ref{thm:VoronoiZone}}
Given $\mathcal{X}$ locally finite, $x \in \mathcal{X}$, an
admissible $A \subset\mathbb{R}^d$,
and $A_0 \subset{\partial A}$, put $\xi^{(5,\ell)}(x, \mathcal{X},
A_0)$ to be the
number of $\ell$-dimensional faces of $v(x)$ if $v(x) \cap A_0 \neq
\varnothing$ and zero otherwise. Define $\xi^{(5,\ell)}(x, \mathcal
{X}, \mathbb{R}
^{d-1})$ similarly.
Now, put
\[
\xi^{(5)}(x, \mathcal{X}, A_0):= \sum
_{l = 0}^{d-1} \xi^{(5,\ell
)}(x, \mathcal{X},
A_0)
\]
and notice that
\[
\operatorname{Co}_\lambda (A_0) = \sum_{x \in{\mathcal{P}}_\lambda }
\xi ^{(5)}(x, {\mathcal{P}}_\lambda, A_0).
\]
We shall show that $\xi^{(5,\ell)}$ satisfies the hypotheses of
Theorems \ref{eq:WLLN} and \ref{eq:Var}, and thus deduce Theorem~\ref
{thm:VoronoiZone}
from \eqref{eq: newexpect} and \eqref{eq: newvar} with $\xi$ set to
$\xi^{(5)}$ and $\gamma$ set to zero (notice that $\xi^{(5)}$ is
homogeneous of order $0$).
The score function $\xi^{(5)}$ is homogeneously stabilizing as at
\eqref{eq:hom}, exponentially stabilizing as at \eqref{eq:expo}, and
satisfies the moment condition \eqref{eq:mom} for all $p \geq1$. This
is because each $\xi^{(5,\ell)}$ with $\ell\in\{0,\ldots,d-1 \}$
has this property. Also, since each $\xi^{(5,\ell)}$ is well
approximated by ${\mathcal{P}}_\lambda $ input on half-spaces for
each $\ell\in\{
0,\ldots,d-1 \}$, it follows that $\xi^{(5)}$ enjoys this property as
well. Thus $\xi^{(5)}$ satisfies the hypotheses of Theorems \ref
{eq:WLLN} and \ref{eq:Var}, concluding the proof of Theorem~\ref
{thm:VoronoiZone}.
\end{pf*}

%\CT[You pointed out a little error below, which is fixed now (a +1 was
%missing). I also modified the conditioning, which is now on $\mbox{PV}_
%\la$ and not %just its volume, which would in fact not be sufficient.]
\begin{pf*}{Proof of Theorem~\ref{thm:IteratedVoronoi}}
We start with the iterated volume $V_\lambda ^{(n)}$. Conditioned on
$\operatorname{PV}_\lambda ^{(1)}$
the first asymptotic equivalence of Theorem~\ref{thm:VoronoiV+S} yields
\[
\mathbb{E} \bigl[V_\lambda ^{(2)}-V_\lambda
^{(1)}\vert \operatorname{PV}_\lambda ^{(1)} \bigr]\sim
c_1 \lambda ^{-{1/d}} {\mathcal{H}}^{d-1} \bigl(
\partial \bigl(\operatorname{PV}_\lambda ^{(1)} \bigr) \bigr).
\]
Taking expectations
and recalling the equivalence $\mathbb{E}S_\lambda (A) \sim c_2
{\mathcal{H}}
^{d-1}(\partial A)$, we obtain
\[
\mathbb{E}V_\lambda ^{(2)} - V(A) \sim c_1 \lambda
^{-{1/d}} (c_2+1) {\mathcal{H}} ^{d-1}(\partial A).
\]
Next,
\begin{eqnarray*}
\mathbb{E}V_\lambda ^{(3)} - V(A) &=& \mathbb{E}\mathbb {E}
\bigl[V_\lambda ^{(3)}-V_\lambda ^{(2)}\vert \operatorname{PV}_\lambda ^{(2)} \bigr]+\mathbb{E}\mathbb{E}
\bigl[V_\lambda ^{(2)}-V_\lambda ^{(1)}\vert \operatorname{PV}_\lambda ^{(1)} \bigr]+\mathbb{E}V_\lambda
^{(1)}-V(A)
\\
&\sim &c_1\lambda ^{-{1/d}} c_{2}^2{
\mathcal{H}}^{d-1}(\partial A)+c_1 \lambda ^{-{1/d}}
c_2 {\mathcal{H}}^{d-1}(\partial A)+c_1 \lambda
^{-{1/d}}{\mathcal{H}} ^{d-1}(\partial A)
\\
&=&c_1 c_{2,2}\lambda ^{-{1/d}}{\mathcal{H}}^{d-1}(
\partial A).
\end{eqnarray*}
Recursively continuing this way proves the desired claim, namely
\[
\mathbb{E}V_\lambda ^{(n)}-V(A)\sim c_{1}
c_{2,n} \lambda ^{-{1/d}} {\mathcal{H}} ^{d-1}({\partial
A}).
\]
The asymptotic equivalences for $\mathbb{E}S_\lambda ^{(n)}$,
$\mathbb{E}H_\lambda ^{\ell,(n)}$ and $\mathbb{E}f_\lambda ^{\ell,(n)}$ follow similarly.
\end{pf*}

%s5
\section{Variance lower bounds}\label{sec:VLB}

%\CT[First paragraph modified, saying that we only show some of the
%results. Is it OK? Also slight modifications below.]
%\CT[This section is still full of $Q$'s. What shall we do with them?]

We complete the proofs of Theorems \ref{thm:VoronoiV+S}--\ref
{thm:VoronoiZone} by proving positivity of the constants appearing in
the variance expressions.
%under the condition that the boundary of the set $A$ contains a subset
%$\Gamma$ of differentiability class $C^2$, which has positive
%$(d-1)$-dimensional Hausdorff measure.
The assumption that $\partial A$ contains a $C^2$-smooth subset with
positive $(d-1)$-dimensional Hausdorff measure is essential for our
following arguments, but we conjecture that this condition can be
relaxed. For example, in \cite{SchulteVoronoi} the author establishes
upper and lower bounds on $\operatorname{Var}[V_\lambda(A)]$ for any
compact convex
set $A$ having non-empty interior, without additional smoothness
assumptions. However, it is unclear (to us) whether the methods of
\cite{SchulteVoronoi} extend to the
more general class of admissible sets $\mathbf{A}(d)$ as well as to
the other
Poisson--Voronoi statistics considered in Theorems \ref
{thm:VoronoiV+S}--\ref{thm:VoronoiZone}.
%and to our more general class of admissible sets $A\in\bA(d)$.

%To simplify the presentation, we show positivity of the constants
%$c_3$ and $c_4$ in Theorem~\ref{thm:VoronoiV+S} only. Positivity of
%$c_6$, $c_8$ and $c_{10}$ can be proved by the %same method.

In what follows, we use the standard Landau notation. More precisely,
for two functions $f,g:[0,\infty)\to\mathbb{R}$ we write
\begin{itemize}
\item $f=o(g)$ if for all $c\in(0,\infty)$ there exists $\lambda
_0>0$ such that for all $\lambda\geq\lambda_0$, $\llvert  f(\lambda)\rrvert   \leq
c\llvert  g(\lambda)\rrvert   $,
\item $f=O(g)$ if there exists $c\in(0,\infty)$ and $\lambda
_0>0$ such that for all $\lambda\geq\lambda_0$, $\llvert  f(\lambda)\rrvert   \leq
c\llvert  g(\lambda)\rrvert   $, and
\item $f=\Omega(g)$ if there exists $c\in(0,\infty)$ and
$\lambda_0>0$ such that for all $\lambda\geq\lambda_0$, $\llvert  f(\lambda)\rrvert   \geq c g(\lambda)$.
\end{itemize}

\subsection*{Positivity of $c_3$ and $c_4$}
Positivity of $c_3$ is
shown in Theorem~2.3 of \cite{Yu} and it remains to consider $c_4$.
For this, recall that $\Gamma\subset\partial A$ is $C^2$-smooth, with
$\mathcal{H}^{d-1}(\Gamma) \in(0, \infty)$. Recalling $A \subset Q$,
subdivide $Q$ into cubes
of edge length $l(\lambda ):= (\lfloor\lambda ^{1/d} \rfloor
)^{-1}$. The
number $L(\lambda )$ of cubes having non-empty intersection with
$\Gamma$ satisfies $L(\lambda ) = \Omega(\lambda ^{(d-1)/d} )$, as otherwise
the cubes would partition $\Gamma$ into $o(\lambda ^{(d-1)/d} )$ sets,
each of
$\mathcal{H}^{d-1}$-measure $O( (\lambda ^{-1/d})^{d-1}) $, which
when $\lambda \to
\infty$ gives
$\mathcal{H}^{d-1}(\Gamma) = 0$, a contradiction.

We find a sub-collection $Q_1,\ldots,Q_M$ of the $L(\lambda )$ cubes
such that
$d(Q_i, Q_j) \geq2 \sqrt{d}  l(\lambda )$ for all $i, j \leq M$,
and $M =
\Omega(\lambda ^{(d-1)/d})$, where $d(Q_i,Q_j)$ stands for the distance
between $Q_i$ and~$Q_j$.\vspace*{1pt} Rotating and translating $Q_i, 1 \leq i \leq
M$, by a distance at most $(\sqrt{d}/2)  l(\lambda )$, if necessary, we
obtain a collection $\widetilde{Q}_1,\ldots,\widetilde{Q}_M$ of
disjoint cubes (with faces
not necessarily parallel to a coordinate plane) such that
\begin{itemize}

\item$d(\widetilde{Q}_i, \widetilde{Q}_j) \geq\sqrt{d}
l(\lambda )$ for all $i, j \leq M$,\vspace*{1pt}

\item$\Gamma$ contains the centre of each $\widetilde{Q}_i$, here denoted
$x_i, 1 \leq i \leq M$.
\end{itemize}
By the assumed differentiability of $\Gamma$, $\Gamma\cap\widetilde
{Q}_i$ is
well approximated
locally around each $x_i$ by the hyperplane $T_i:= T_{x_i}$ tangent
to $\Gamma$ at $x_i$. By the $C^2$-assumption, the approximation is
uniform over all $1 \leq i \leq M$. % (this assumption seems to be
%missing in \cite{Yu}).
Making a further rotation of $\widetilde{Q}_i$, if necessary, we may
assume that
$T_i$ partitions $\widetilde{Q}_i$ into
congruent rectangular solids. Let $T_i$ coincide with the hyperplane
$\mathbb{R}^{d-1}$. Without loss of generality, we assume $\partial A
\subset
\mathbb{R}^{d-1} \times(-\infty, 0]$, that is,
$\partial A$ is `beneath' $T_i$.

We now exhibit a configuration of Poisson
points ${\mathcal{P}}_\lambda $ which has strictly positive
probability and for which
$S_\lambda (A)$ has variability bounded below by $\Omega( \lambda
^{-(d-1)/d}\mathcal{H}^{d-1}(\Gamma))$. Let $\epsilon:= \epsilon
(\lambda ):=
l(\lambda )/28$ and sub-divide each
$\widetilde{Q}_i, 1 \leq i \leq M$, into $28^d$ sub-cubes of edge length
$\epsilon$.
Sub-cubes within Hausdorff distance $4 \epsilon$ of $\partial
\widetilde{Q}_i$
are called `boundary' sub-cubes; if a sub-cube is not a boundary sub-cube
then we call it an interior sub-cube. If each boundary sub-cube in
$\widetilde{Q}
_i$ contains a point from ${\mathcal{P}}_\lambda $, then the geometry
of the Voronoi
cells with centres in $\widetilde{Q}_i$ and
distant more than $4 \epsilon$ from $\partial\widetilde{Q}_i$ is not altered
by point configurations outside
$\widetilde{Q}_i$ (see, e.g., \cite{Penrose07b}).

We assume that $x_i$ coincides with the origin and we recall that
$\partial A\subset\mathbb{R}^{d-1}\times(-\infty,0]$ so that points near
$\partial A$
may be parametrized by a pair in $\mathbb{R}^{d-1}\times(-\infty,0]$). By
$2(\mathbb{Z}^{d-1})$ we mean the set of all points in $\mathbb
{R}^{d-1}$ having
integer coordinates of even parity.
Consider the sub-cubes $\widetilde{Q}_i$ having the following properties:
\begin{enumerate}[(b')]

\item[(a)] the boundary sub-cubes each contain at least one point from
${\mathcal{P}}_\lambda $,

\item[(b)] ${\mathcal{P}}_\lambda \cap B_{\epsilon/100 }(( \epsilon
j, \pm\epsilon
))$ consists of a singleton
for $j \in2(\mathbb{Z}^{d-1}),\llvert  j\rrvert    \leq10$, %$j \in\{0, \pm2, \pm4,...,
%\pm10 \}$,
or
\item[(b$'$)] ${\mathcal{P}}_\lambda \cap B_{\epsilon/100 }(( \epsilon
j, \epsilon
/100))$ consists of a singleton
for $j \in2(\mathbb{Z}^{d-1}),\llvert  j\rrvert    \leq10$
%$j \in\{0, \pm2, \pm4,... \pm10 \}$
and also ${\mathcal{P}}_\lambda \cap B_{\epsilon/100 }(( \epsilon j,
-\epsilon
/100))$ consists of a singleton
for $j=0$ and $j \in2(\mathbb{Z}^{d-1})+1,\llvert  j\rrvert    \leq10$,
%$j \in\{0, \pm1, \pm3,...,\pm9 \},$

\item[(c)] ${\mathcal{P}}_\lambda $ puts no other points in
$\widetilde{Q}_i$.
\end{enumerate}
(We remark that the choice of the constants $28$ and $100$ is arbitrary
and that we could have used any sufficiently large number.) Events (b)
and (b$'$), which each involve $22$ singletons in $22$ small balls,
happen with the same probability, which is small but bounded away from
zero uniformly in $\lambda $, since $\kappa\equiv1$.

Re-labelling if necessary, let $I:= \{1,\ldots,K\}$ be the indices of
cubes $\widetilde{Q}_i$ having properties (a)--(c). It is easily
checked that
the probability a given $\widetilde{Q}_i, 1 \leq i \leq M$, satisfies property
(a) is strictly
positive, uniformly in $\lambda $. This is also true for properties
(b)--(c), showing that
\begin{equation}
\label{EK2} \mathbb{E}K = \Omega \bigl(\lambda ^{(d-1)/d} \bigr).
\end{equation}
Abusing notation, let $\mathcal{Q}:= \bigcup_{i=1}^K \widetilde
{Q}_{i}$ and put $\mathcal{Q}^c:= [0,1]^d \setminus\mathcal{Q}$. Let $\mathcal{F}_\lambda $ be the
$\sigma$-algebra
determined by the random set $I$, the positions of points of ${\mathcal
{P}}_\lambda $
in all boundary sub-cubes,
and the positions of points ${\mathcal{P}}_\lambda $ in $\mathcal
{Q}^c$. Let $U_i, 1 \leq i
\leq M$, be the union of the interior sub-cubes
in $\widetilde{Q}_i$. If $d=2$, we notice that if (b) happens, then
the surface
$\partial A_\lambda \cap U_i$ contains nearly horizontal edges and the
total length
of these edges is generously bounded above by $30 \epsilon$. Indeed,
if (b) happens, the $11$
cells centered at the points in
${\mathcal{P}}_\lambda \cap B_{\epsilon/100 }(( \epsilon j, -
\epsilon))$, $j \in
\{0, \pm2, \pm4,\ldots, \pm10 \}$,
contribute to $\partial(\operatorname{PV}_\lambda (A))$ a length roughly bounded by
the width
of $U_i$ plus some negligible corrections.
On the other hand,
if (b$'$) happens then $\partial A_\lambda \cap U_i$ contains $10$ sharp
peaks, with abscissas roughly equal to $\{\pm1, \pm3,\ldots, \pm9 \}
$. In fact, it is easily checked that
$\partial A_\lambda \cap U_i$ contains at least $18$ `long', nearly
vertical edges of length at least
$2 \epsilon$, giving a total edge length of at least $36 \epsilon$.
A similar situation holds in higher dimensions $d\geq3$. % with $
%\epsilon$ replaced by $\epsilon^{d-1}$.

Conditional on $\mathcal{F}_\lambda $, $\mathcal{H}^{d-1}(\partial
A_\lambda \cap\widetilde{Q}_i)$ has
variability $\Omega( \epsilon^{2(d-1)}) = \Omega( \lambda ^{-2 + 2/d}
)$, uniformly in $i \in I$,
that is,
\begin{equation}
\label{VV2} \operatorname{Var} \bigl[\mathcal {H}^{d-1}(\partial
A_\lambda \cap\widetilde{Q}_i) \vert \mathcal{F}_\lambda
\bigr] = \Omega \bigl( \lambda ^{-2 + 2/d} \bigr), \qquad i \in I.
\end{equation}
By the conditional
variance formula,
\begin{eqnarray*}
\operatorname{Var} \bigl[S_\lambda (A) \bigr] &=& \operatorname{Var}
\bigl[ \mathbb{E}\bigl[S_\lambda (A) \vert \mathcal{F}_\lambda \bigr]\bigr] +
\mathbb{E} \bigl[ \operatorname{Var} \bigl[ S_\lambda (A) \vert \mathcal
{F}_\lambda \bigr] \bigr]
\\
&\geq&\mathbb{E} \bigl[ \operatorname{Var} \bigl[ S_\lambda (A)\vert
\mathcal {F}_\lambda \bigr] \bigr]
\\
&=&\mathbb{E} \bigl[\operatorname{Var} \bigl[\mathcal{H}^{d-1}(\partial
A_\lambda \cap\mathcal{Q}) + \mathcal{H}^{d-1} \bigl(\partial
A_\lambda \cap\mathcal{Q}^c \bigr) \vert
\mathcal{F}_\lambda \bigr] \bigr].
\end{eqnarray*}
Given $\mathcal{F}_\lambda $, the Poisson--Voronoi mosaic of
${\mathcal{P}}_\lambda $ admits
variability only inside $\mathcal{Q}$, that is to say, given $\mathcal
{F}_\lambda $, we
have $\mathcal{H}^{d-1}(\partial A_\lambda \cap\mathcal{Q}^c)$ is
constant. Thus,
\begin{eqnarray}
\operatorname{Var} \bigl[S_\lambda (A) \bigr] &\geq&\mathbb{E} \bigl[
\operatorname {Var} \bigl[ \mathcal{H}^{d-1}(\partial A_\lambda
\cap\mathcal{Q}) \vert \mathcal{F}_\lambda \bigr] \bigr]
\nonumber
\\
&=& \mathbb{E} \biggl[\operatorname{Var} \biggl[ \sum
_{i \in I} \mathcal {H}^{d-1}(\partial A_\lambda
\cap\widetilde{Q}_i) \Big| \mathcal{F}_\lambda \biggr]
\biggr]
\label{VV2a}
\\
& =& \mathbb{E}\sum_{i \in I} \operatorname{Var} \bigl[
\mathcal {H}^{d-1}(\partial A_\lambda \cap\widetilde{Q}_i)
\vert \mathcal{F} _\lambda \bigr],\nonumber
\end{eqnarray}
since, given $\mathcal{F}_\lambda $, $\mathcal{H}^{d-1}(\partial
A_\lambda \cap\widetilde{Q}_i), i \in
I$, are
independent. By \eqref{EK2} and \eqref{VV2}, we have
\[
\operatorname{Var} \bigl[ S_\lambda (A) \bigr] \geq c \lambda
^{-2 + 2/d} \mathbb{E}[K] = \Omega \bigl(\lambda ^{-(d -
1)/d} \bigr)
\]
with some finite constant $c \in(0, \infty)$, concluding the proof
that $c_4$ is positive.

\subsection*{Positivity of $c_6$ and $c_8$}

The general idea is to show
that configuration (b$'$) generates a surface which has more variability
(both in terms of complexity and measure) than the surface generated by
configuration (b). The details go as follows.
For $\ell\in\{0,1,\ldots,d-1 \}$ and $i \in I$, put $S_{\ell, i}:=
(\operatorname{skel}_\ell(\operatorname{PV}_\lambda (A))) \cap U_i$, noting that
$\partial(\operatorname{PV}_\lambda (A)) \cap U_i = S_{d-1, i}$ (recall the notation
introduced in the discussion around equation \eqref{EK2}). Let
$S_{\ell, i}(b)$ be the $\ell$-dimensional skeleton arising from
configuration (b) and define $S_{\ell, i}(b')$ similarly. Henceforth,
without loss of generality we fix $i = 1$ and write $S_\ell$ for
$S_{\ell,1}$. Observe that $S_{d-1}(b)$
consists of a single $(d-1)$-dimensional facet $f$ which is nearly a
hypercube of dimension $d-1$ (and nearly a horizontal edge when $d = 2$).
Also, $S_{d-2}(b)$ is the union of $2(d-1)$ faces, each of which is
nearly a hypercube of dimension $d-2$.

On the other hand, $S_{d-1}(b')$ is the union of $(d-1)$-dimensional
facets $F_j, 1 \leq j \leq2(d-1)$, whose union forms the boundary of a solid
hyper-pyramid in $\mathbb{R}^d$ whose
base is a translate, up to a negligible perturbation, of $S_{d-1}(b)$.
The boundary of the surface $\bigcup_{j = 1}^{2(d-1)} F_j$ is of
dimension $d-2$
and is the union of $2(d-1)$ faces, each of which is nearly a hypercube
of dimension $d-2$.
In fact, the boundary of the surface $\bigcup_{j = 1}^{2(d-1)} F_j$ is a
translate,\vspace*{-1pt} also up to a negligible perturbation, of $S_{d-2}(b)$; we
thus denote the boundary of $\bigcup_{j = 1}^{2(d-1)} F_j$ by
$\tilde{S}_{d-2}(b')$. In other words, we have that
\[
S_{d-1}(b) = S_{d-2}(b) \cup(\operatorname{int} f)
\]
and
\[
S_{d-1} \bigl(b' \bigr) = \tilde{S}_{d-2}
\bigl(b' \bigr) \cup \Biggl( \bigcup_{j =
1}^{2(d-1)}
F_j \setminus\tilde{S}_{d-2}(b) \Biggr).
\]
Now, $S_{d-2}(b)$ and $\tilde{S}_{d-2}(b')$ are indistinguishable from
the viewpoint of their combinatorial complexity, as measured
by their lower-dimensional skeletons. Moreover, they are nearly
indistinguishable from a measure theoretic point of view, since the
$\mathcal{H}
^{\ell}$-measure
of their $\ell$-skeletons nearly coincide (modulo negligible
corrections). On the other hand, the open facet $\operatorname{int} f$ differs
significantly from\vspace*{-1pt}
$\bigcup_{j = 1}^{2(d-1)} F_j \setminus\tilde{S}_{d-2}(b)$ in terms of
both combinatorial complexity and measure. Indeed,
the $\mathcal{H}^{\ell}$-measure of
the $\ell$-skeleton of the latter (facets of a pyramid) is strictly
larger than the $\mathcal{H}^{\ell}$-measure of $\operatorname{int} f$ (the
base of
the pyramid). Likewise, for $\ell\in\{0,1,\ldots,d-2 \}$, the single
facet $\operatorname{int} f$ has no $\ell$-dimensional faces, whereas $\bigcup
_{j = 1}^{2(d-1)} F_j \setminus\tilde{S}_{d-2}(b)$
has a non-zero number of $\ell$-dimensional faces. These arguments
apply to all skeletons $S_{d-1, i}(b), i \in I$. By following nearly
verbatim the arguments
showing that $c_4$ is positive, we get that $c_6$ and $c_8$ are positive.

\subsection*{Positivity of $c_{10}$}

We have that $\operatorname{Co}_\lambda (A_0)$
is defined in terms of $f_\lambda ^{(\ell)}(A), \ell\in\{0,1,\ldots,d-1
\}$,
and it suffices to note that configuration (b$'$) leads to a complexity
which is strictly larger than the complexity arising from
configuration (b). We now follow the arguments that $c_4$ is strictly positive.

\section*{Acknowledgements}
We thank two anonymous referees for comments leading to an improved
exposition. Christoph Th\"ale has been supported
by the German Research Foundation (DFG) via SFB-TR 12 whereas Joseph
Yukich has been supported
in part by NSF grant DMS-1406410. Joseph Yukich gratefully thanks
the Faculty of Mathematics at Ruhr University Bochum for its kind
hospitality and support.

%\begin{appendix}
%\section{}
%\end{appendix}

%\section*{Acknowledgements}

%\begin{supplement}%[id=suppA]
%\slink[doi]{10.3150/00-BEJXXXXSUPP} %[doi,text={...}] - jei reikia
%suskaldyti doi
%\sfilename{BEJ000\_supp.pdf}
%\end{supplement}

% imsref loaded by imikolaityte, 2015-08-17 08:29:06
%

\printhistory

\begin{thebibliography}{31}
% pybtex-1.42. Style name=bej, version=1.42, label_style=nolabel,
%sorting_style=complex, cfg=None, language=None.

%b1 ###
\bibitem{BaiHwangEtal01}
%
\begin{barticle}[mr]
\bauthor{\bsnm{Bai}, \bfnm{Zhi-Dong}\binits{Z.-D.}},
\bauthor{\bsnm{Hwang}, \bfnm{Hsien-Kuei}\binits{H.-K.}},
\bauthor{\bsnm{Liang}, \bfnm{Wen-Qi}\binits{W.-Q.}} \AND
\bauthor{\bsnm{Tsai}, \bfnm{Tsung-Hsi}\binits{T.-H.}}
(\byear{2001}).
\btitle{Limit theorems for the number of maxima in random samples from
planar regions}.
\bjournal{Electron. J. Probab.}
\bvolume{6}
\bpages{no. 3, 41 pp. (electronic)}.
\bid{doi={10.1214/EJP.v6-76}, issn={1083-6489}, mr={1816046}}
\bptnote{check pages}%
\end{barticle}
%
\iffalse\OrigBibText
Bai, Z.-D., H.-K. Hwang, W.-Q. Liang and Tsai, T.-H. (2001) Limit
theorems for the number of maxima in random samples from planar
regions. \textit{Electron. J. Probab.} \textbf{6}, 1--41.
\endOrigBibText\fi
\bptok{imsref}%
% NOT OUTPUTTED:
% doi = http://dx.doi.org/10.1214/EJP.v6-76
% fjournal = Electronic Journal of Probability
\endbibitem

%b2 ###
\bibitem{BarbourXia01}
%
\begin{barticle}[mr]
\bauthor{\bsnm{Barbour}, \bfnm{A. D.}\binits{A.D.}} \AND
\bauthor{\bsnm{Xia}, \bfnm{A.}\binits{A.}}
(\byear{2001}).
\btitle{The number of two-dimensional maxima}.
\bjournal{Adv. in Appl. Probab.}
\bvolume{33}
\bpages{727--750}.
\bid{doi={10.1239/aap/1011994025}, issn={0001-8678}, mr={1875775}}
\end{barticle}
%
\iffalse\OrigBibText
Barbour, A.D. and Xia, A. (2001) The number of two dimensional maxima.
\textit{Adv. in Appl. Probab.} \textbf{33}, 727--750.
\endOrigBibText\fi
\bptok{imsref}%
% NOT OUTPUTTED:
% number = 4
% doi = http://dx.doi.org/10.1239/aap/1011994025
% coden = AAPBBD
% fjournal = Advances in Applied Probability
\endbibitem

%b3 ###
\bibitem{BarbourXia06}
%
\begin{barticle}[mr]
\bauthor{\bsnm{Barbour}, \bfnm{A. D.}\binits{A.D.}} \AND
\bauthor{\bsnm{Xia}, \bfnm{Aihua}\binits{A.}}
(\byear{2006}).
\btitle{Normal approximation for random sums}.
\bjournal{Adv. in Appl. Probab.}
\bvolume{38}
\bpages{693--728}.
\bid{doi={10.1239/aap/1158684998}, issn={0001-8678}, mr={2256874}}
\end{barticle}
%
\iffalse\OrigBibText
Barbour, A.D. and Xia, A. (2006) Normal approximation for random sums.
\textit{Adv. in Appl. Probab.} \textbf{38}, 693--728.
\endOrigBibText\fi
\bptok{imsref}%
% NOT OUTPUTTED:
% number = 3
% doi = http://dx.doi.org/10.1239/aap/1158684998
% coden = AAPBBD
% fjournal = Advances in Applied Probability
\endbibitem

%b4 ###
\bibitem{BaryshnikovYukich2005}
%
\begin{barticle}[mr]
\bauthor{\bsnm{Baryshnikov}, \bfnm{Yu}\binits{Y.}} \AND
\bauthor{\bsnm{Yukich}, \bfnm{J. E.}\binits{J.E.}}
(\byear{2005}).
\btitle{Gaussian limits for random measures in geometric probability}.
\bjournal{Ann. Appl. Probab.}
\bvolume{15}
\bpages{213--253}.
\bid{doi={10.1214/105051604000000594}, issn={1050-5164}, mr={2115042}}
\end{barticle}
%
\iffalse\OrigBibText
Baryshnikov, Y. and Yukich, J.E. (2005) Gaussian limits for random
measures in geometric probability. \textit{Ann. Appl. Probab.} \textbf
{15}, 213--253.
\endOrigBibText\fi
\bptok{imsref}%
% NOT OUTPUTTED:
% number = 1A
% doi = http://dx.doi.org/10.1214/105051604000000594
% fjournal = The Annals of Applied Probability
\endbibitem

%b5 ###
\bibitem{Devroye}
%
\begin{barticle}[mr]
\bauthor{\bsnm{Devroye}, \bfnm{Luc}\binits{L.}}
(\byear{1993}).
\btitle{Records, the maximal layer, and uniform distributions in
monotone sets}.
\bjournal{Comput. Math. Appl.}
\bvolume{25}
\bpages{19--31}.
\bid{doi={10.1016/0898-1221(93)90195-2}, issn={0898-1221}, mr={1199909}}
\end{barticle}
%
\iffalse\OrigBibText
Devroye, L. (1993): Records, the maximal layer, and uniform
distributions in monotone sets. \textit{Computers Math. Appl.} \textbf
{25}, 19--31.
\endOrigBibText\fi
\bptok{imsref}%
% NOT OUTPUTTED:
% number = 5
% doi = http://dx.doi.org/10.1016/0898-1221(93)90195-2
% coden = CMAPDK
% fjournal = Computers \& Mathematics with Applications. An
%International Journal
\endbibitem

%b6 ###
\bibitem{FedererPaper}
%
\begin{barticle}[mr]
\bauthor{\bsnm{Federer}, \bfnm{Herbert}\binits{H.}}
(\byear{1959}).
\btitle{Curvature measures}.
\bjournal{Trans. Amer. Math. Soc.}
\bvolume{93}
\bpages{418--491}.
\bid{issn={0002-9947}, mr={0110078}}
\end{barticle}
%
\iffalse\OrigBibText
Federer, H. (1959) Curvature measures. \textit{Trans. Am. Math. Soc.}
\textbf{93}, 418--491.
\endOrigBibText\fi
\bptok{imsref}%
% NOT OUTPUTTED:
% fjournal = Transactions of the American Mathematical Society
\endbibitem

%b7 ###
\bibitem{Federer}
%
\begin{bbook}[mr]
\bauthor{\bsnm{Federer}, \bfnm{Herbert}\binits{H.}}
(\byear{1969}).
\btitle{Geometric Measure Theory}.
\bseries{Die Grundlehren der Mathematischen Wissenschaften}
\bvolume{153}.
\blocation{New York}:
\bpublisher{Springer}.
\bid{mr={0257325}}
\end{bbook}
%
\iffalse\OrigBibText
Federer, H. (1969) \textit{Geometric Measure Theory}. Springer, Berlin.
\endOrigBibText\fi
\bptok{imsref}%
% NOT OUTPUTTED:
% fpage = xiv+676
\endbibitem

%b8 ###
\bibitem{Fremlin}
%
\begin{barticle}[mr]
\bauthor{\bsnm{Fremlin}, \bfnm{D. H.}\binits{D.H.}}
(\byear{1997}).
\btitle{Skeletons and central sets}.
\bjournal{Proc. Lond. Math. Soc. (3)}
\bvolume{74}
\bpages{701--720}.
\bid{doi={10.1112/S0024611597000233}, issn={0024-6115}, mr={1434446}}
\end{barticle}
%
\iffalse\OrigBibText
Fremlin, D.H. (1997) Skeletons and central sets. \textit{Proc. London
Math. Soc.} \textbf{74}, 701--720.
\endOrigBibText\fi
\bptok{imsref}%
% NOT OUTPUTTED:
% number = 3
% doi = http://dx.doi.org/10.1112/S0024611597000233
% coden = PLMTAL
% fjournal = Proceedings of the London Mathematical Society. Third
%Series
\endbibitem

%b9 ###
\bibitem{HevelingReitzner}
%
\begin{barticle}[mr]
\bauthor{\bsnm{Heveling}, \bfnm{Matthias}\binits{M.}} \AND
\bauthor{\bsnm{Reitzner}, \bfnm{Matthias}\binits{M.}}
(\byear{2009}).
\btitle{Poisson--{V}oronoi approximation}.
\bjournal{Ann. Appl. Probab.}
\bvolume{19}
\bpages{719--736}.
\bid{doi={10.1214/08-AAP561}, issn={1050-5164}, mr={2521886}}
\end{barticle}
%
\iffalse\OrigBibText
Heveling, M. and Reitzner, M. (2009) Poisson-Voronoi approximation.
\textit{Ann. Appl. Probab.} \textbf{19}, 719--736.
\endOrigBibText\fi
\bptok{imsref}%
% NOT OUTPUTTED:
% number = 2
% doi = http://dx.doi.org/10.1214/08-AAP561
% fjournal = The Annals of Applied Probability
\endbibitem

%b10 ###
\bibitem{HLW}
%
\begin{barticle}[mr]
\bauthor{\bsnm{Hug}, \bfnm{Daniel}\binits{D.}},
\bauthor{\bsnm{Last}, \bfnm{G{\"u}nter}\binits{G.}} \AND
\bauthor{\bsnm{Weil}, \bfnm{Wolfgang}\binits{W.}}
(\byear{2004}).
\btitle{A local {S}teiner-type formula for general closed sets and
applications}.
\bjournal{Math. Z.}
\bvolume{246}
\bpages{237--272}.
\bid{doi={10.1007/s00209-003-0597-9}, issn={0025-5874}, mr={2031455}}
\end{barticle}
%
\iffalse\OrigBibText
Hug, D., Last, G. and Weil, W. (2004) A local Steiner-type formula for
general closed sets and applications. \textit{Math. Z.} \textbf{246},
237--272.
\endOrigBibText\fi
\bptok{imsref}%
% NOT OUTPUTTED:
% number = 1-2
% doi = http://dx.doi.org/10.1007/s00209-003-0597-9
% coden = MAZEAX
% fjournal = Mathematische Zeitschrift
\endbibitem

%b11 ###
\bibitem{HwangTsai10}
%
\begin{barticle}[mr]
\bauthor{\bsnm{Hwang}, \bfnm{Hsien-Kuei}\binits{H.-K.}} \AND
\bauthor{\bsnm{Tsai}, \bfnm{Tsung-Hsi}\binits{T.-H.}}
(\byear{2010}).
\btitle{Multivariate records based on dominance}.
\bjournal{Electron. J. Probab.}
\bvolume{15}
\bpages{1863--1892}.
\bid{doi={10.1214/EJP.v15-825}, issn={1083-6489}, mr={2738341}}
\bptnote{check pages}%
\end{barticle}
%
\iffalse\OrigBibText
Hwang, H.-K. and Tsai, T.-H. (2010) Multivariate records based on dominance.
\textit{Electron. J. Probab.} \textbf{15}, 1863--1892.
\endOrigBibText\fi
\bptok{imsref}%
% NOT OUTPUTTED:
% doi = http://dx.doi.org/10.1214/EJP.v15-825
% fjournal = Electronic Journal of Probability
\endbibitem

%b12 ###
\bibitem{KestenLee}
%
\begin{barticle}[mr]
\bauthor{\bsnm{Kesten}, \bfnm{Harry}\binits{H.}} \AND
\bauthor{\bsnm{Lee}, \bfnm{Sungchul}\binits{S.}}
(\byear{1996}).
\btitle{The central limit theorem for weighted minimal spanning trees
on random points}.
\bjournal{Ann. Appl. Probab.}
\bvolume{6}
\bpages{495--527}.
\bid{doi={10.1214/aoap/1034968141}, issn={1050-5164}, mr={1398055}}
\end{barticle}
%
\iffalse\OrigBibText
Kesten, H. and Lee, S. (1996) The central limit theorem for weighted
minimal spanning trees on random points. \textit{Ann. Appl. Probab.}
\textbf{6}, 495--527.
\endOrigBibText\fi
\bptok{imsref}%
% NOT OUTPUTTED:
% number = 2
% doi = http://dx.doi.org/10.1214/aoap/1034968141
% fjournal = The Annals of Applied Probability
\endbibitem

%b13 ###
\bibitem{KT}
%
\begin{barticle}[mr]
\bauthor{\bsnm{Khmaladze}, \bfnm{Estate}\binits{E.}} \AND
\bauthor{\bsnm{Toronjadze}, \bfnm{N.}\binits{N.}}
(\byear{2001}).
\btitle{On the almost sure coverage property of {V}oronoi
tessellation: The {$\Bbb R\sp1$} case}.
\bjournal{Adv. in Appl. Probab.}
\bvolume{33}
\bpages{756--764}.
\bid{doi={10.1239/aap/1011994027}, issn={0001-8678}, mr={1875777}}
\end{barticle}
%
\iffalse\OrigBibText
E. Khmaladze and N. Toronjadze (2001) On the almost sure coverage
property of Voronoi tessellation: the $\mathbb{R}^1$ case. \textit
{Adv. Appl.
Probab.} \textbf{33}, 756--764.
\endOrigBibText\fi
\bptok{imsref}%
% NOT OUTPUTTED:
% number = 4
% doi = http://dx.doi.org/10.1239/aap/1011994027
% coden = AAPBBD
% fjournal = Advances in Applied Probability
\endbibitem

%b14 ###
\bibitem{KiderlenRataj}
%
\begin{barticle}[mr]
\bauthor{\bsnm{Kiderlen}, \bfnm{Markus}\binits{M.}} \AND
\bauthor{\bsnm{Rataj}, \bfnm{Jan}\binits{J.}}
(\byear{2006}).
\btitle{On infinitesimal increase of volumes of morphological transforms}.
\bjournal{Mathematika}
\bvolume{53}
\bpages{103--127}.
\bid{doi={10.1112/S002557930000005X}, issn={0025-5793}, mr={2304055}}
\bptnote{check pages}%
\end{barticle}
%
\iffalse\OrigBibText
Kiderlen, M. and Rataj, J. (2006) On infinitesimal increase of volumes
of morphological
transforms. \textit{Mathematika} \textbf{53}, 103--127.
\endOrigBibText\fi
\bptok{imsref}%
% NOT OUTPUTTED:
% number = 1
% doi = http://dx.doi.org/10.1112/S002557930000005X
% coden = MTKAAB
% fjournal = Mathematika. A Journal of Pure and Applied Mathematics
\endbibitem

%b15 ###
\bibitem{Lee}
%
\begin{barticle}[mr]
\bauthor{\bsnm{Lee}, \bfnm{Sungchul}\binits{S.}}
(\byear{1997}).
\btitle{The central limit theorem for {E}uclidean minimal spanning
trees. {I}}.
\bjournal{Ann. Appl. Probab.}
\bvolume{7}
\bpages{996--1020}.
\bid{doi={10.1214/aoap/1043862422}, issn={1050-5164}, mr={1484795}}
\end{barticle}
%
\iffalse\OrigBibText
Lee, S. (1997) The central limit theorem for Euclidean minimum spanning
trees. \textit{Ann. Appl. Probab.} \textbf{7}, 996--1020.
\endOrigBibText\fi
\bptok{imsref}%
% NOT OUTPUTTED:
% number = 4
% doi = http://dx.doi.org/10.1214/aoap/1043862422
% fjournal = The Annals of Applied Probability
\endbibitem

%b16 ###
\bibitem{Matousek}
%
\begin{bbook}[mr]
\bauthor{\bsnm{Matou{\v{s}}ek}, \bfnm{Ji\v{r}\'{\i}}\binits{J.}}
(\byear{2002}).
\btitle{Lectures on Discrete Geometry}.
\bseries{Graduate Texts in Mathematics}
\bvolume{212}.
\blocation{New York}:
\bpublisher{Springer}.
\bid{doi={10.1007/978-1-4613-0039-7}, mr={1899299}}
\end{bbook}
%
\iffalse\OrigBibText
Matou\v{s}ek, J. (2002) \textit{Lectures on Discrete Geometry}.
Springer, Berlin.
\endOrigBibText\fi
\bptok{imsref}%
% NOT OUTPUTTED:
% doi = http://dx.doi.org/10.1007/978-1-4613-0039-7
% isbn = 0-387-95373-6
% fpage = xvi+481
\endbibitem

%b17 ###
\bibitem{Mattila}
%
\begin{bbook}[mr]
\bauthor{\bsnm{Mattila}, \bfnm{Pertti}\binits{P.}}
(\byear{1995}).
\btitle{Geometry of Sets and Measures in {E}uclidean Spaces: Fractals and Rectifiability}.
\bseries{Cambridge Studies in Advanced Mathematics}
\bvolume{44}.
\blocation{Cambridge}:
\bpublisher{Cambridge Univ. Press}.
\bid{doi={10.1017/CBO9780511623813}, mr={1333890}}
\end{bbook}
%
\iffalse\OrigBibText
Mattila, P. (1995) \textit{Geometry of Sets and Measures in Euclidean
Spaces}. Cambridge University Press, Cambridge.
\endOrigBibText\fi
\bptok{imsref}%
% NOT OUTPUTTED:
% doi = http://dx.doi.org/10.1017/CBO9780511623813
% isbn = 0-521-46576-1; 0-521-65595-1
% fpage = xii+343
\endbibitem

%b18 ###
\bibitem{MY}
%
\begin{barticle}[mr]
\bauthor{\bsnm{McGivney}, \bfnm{K.}\binits{K.}} \AND
\bauthor{\bsnm{Yukich}, \bfnm{J. E.}\binits{J.E.}}
(\byear{1999}).
\btitle{Asymptotics for {V}oronoi tessellations on random samples}.
\bjournal{Stochastic Process. Appl.}
\bvolume{83}
\bpages{273--288}.
\bid{doi={10.1016/S0304-4149(99)00035-6}, issn={0304-4149}, mr={1708209}}
\end{barticle}
%
\iffalse\OrigBibText
McGivney, K. and Yukich, J.E. (1999) Asymptotics for Voronoi
tessellations. \textit{Stochastic Process. Appl.} \textbf{83}, 273--288.
\endOrigBibText\fi
\bptok{imsref}%
% NOT OUTPUTTED:
% number = 2
% doi = http://dx.doi.org/10.1016/S0304-4149(99)00035-6
% coden = STOPB7
% fjournal = Stochastic Processes and their Applications
\endbibitem

%b19 ###
\bibitem{Mo}
%
\begin{bbook}[auto:parserefs-M02]
\bauthor{\bsnm{Molchanov}, \bfnm{I.}\binits{I.}}
(\byear{1997}).
\btitle{Statistics of the Boolean Model for Practitioners and Mathematicians}.
\blocation{Chichester}:
\bpublisher{Wiley}.
\end{bbook}
%
\iffalse\OrigBibText
Molchanov, I. (1997) \textit{Statistics of the Boolean model for
Practitioners and Mathematicians}. John Wiley and Sons, Chichester.
\endOrigBibText\fi
\bptok{imsref}%
\endbibitem

%b20 ###
\bibitem{Penrose07a}
%
\begin{barticle}[mr]
\bauthor{\bsnm{Penrose}, \bfnm{Mathew D.}\binits{M.D.}}
(\byear{2007}).
\btitle{Gaussian limits for random geometric measures}.
\bjournal{Electron. J. Probab.}
\bvolume{12}
\bpages{989--1035 (electronic)}.
\bid{doi={10.1214/EJP.v12-429}, issn={1083-6489}, mr={2336596}}
\bptnote{check pages}%
\end{barticle}
%
\iffalse\OrigBibText
Penrose, M.D. (2007) Gaussian limits for random geometric measures.
\textit{Electron. J. Probab.} \textbf{12}, 989--1035.
\endOrigBibText\fi
\bptok{imsref}%
% NOT OUTPUTTED:
% doi = http://dx.doi.org/10.1214/EJP.v12-429
% fjournal = Electronic Journal of Probability
\endbibitem

%b21 ###
\bibitem{Penrose07b}
%
\begin{barticle}[mr]
\bauthor{\bsnm{Penrose}, \bfnm{Mathew D.}\binits{M.D.}}
(\byear{2007}).
\btitle{Laws of large numbers in stochastic geometry with statistical
applications}.
\bjournal{Bernoulli}
\bvolume{13}
\bpages{1124--1150}.
\bid{doi={10.3150/07-BEJ5167}, issn={1350-7265}, mr={2364229}}
\end{barticle}
%
\iffalse\OrigBibText
Penrose, M.D. (2007) Laws of large numbers in stochastic geometry with
statistical applications. \textit{Bernoulli} \textbf{13}, 1124--1150.
\endOrigBibText\fi
\bptok{imsref}%
% NOT OUTPUTTED:
% number = 4
% doi = http://dx.doi.org/10.3150/07-BEJ5167
% fjournal = Bernoulli. Official Journal of the Bernoulli Society for
%Mathematical Statistics and Probability
\endbibitem

%b22 ###
\bibitem{PenroseYukich01}
%
\begin{barticle}[mr]
\bauthor{\bsnm{Penrose}, \bfnm{Mathew D.}\binits{M.D.}} \AND
\bauthor{\bsnm{Yukich}, \bfnm{J. E.}\binits{J.E.}}
(\byear{2001}).
\btitle{Central limit theorems for some graphs in computational geometry}.
\bjournal{Ann. Appl. Probab.}
\bvolume{11}
\bpages{1005--1041}.
\bid{doi={10.1214/aoap/1015345393}, issn={1050-5164}, mr={1878288}}
\end{barticle}
%
\iffalse\OrigBibText
Penrose, M.D. and Yukich, J.E. (2001) Central limit theorems for some
graphs in computational geometry. \textit{Ann. Appl. Probab.} \textbf
{11}, 1005--1041.
\endOrigBibText\fi
\bptok{imsref}%
% NOT OUTPUTTED:
% number = 4
% doi = http://dx.doi.org/10.1214/aoap/1015345393
% fjournal = The Annals of Applied Probability
\endbibitem

%b23 ###
\bibitem{PenroseYukich03}
%
\begin{barticle}[mr]
\bauthor{\bsnm{Penrose}, \bfnm{Mathew D.}\binits{M.D.}} \AND
\bauthor{\bsnm{Yukich}, \bfnm{J. E.}\binits{J.E.}}
(\byear{2003}).
\btitle{Weak laws of large numbers in geometric probability}.
\bjournal{Ann. Appl. Probab.}
\bvolume{13}
\bpages{277--303}.
\bid{doi={10.1214/aoap/1042765669}, issn={1050-5164}, mr={1952000}}
\end{barticle}
%
\iffalse\OrigBibText
Penrose, M.D. and Yukich, J.E. (2003) Weak laws of large numbers in
geometric probability. \textit{Ann. Appl. Probab.} \textbf{13}, 277--303.
\endOrigBibText\fi
\bptok{imsref}%
% NOT OUTPUTTED:
% number = 1
% doi = http://dx.doi.org/10.1214/aoap/1042765669
% fjournal = The Annals of Applied Probability
\endbibitem

%b24 ###
\bibitem{PenroseYukich05}
%
\begin{bincollection}[mr]
\bauthor{\bsnm{Penrose}, \bfnm{Mathew D.}\binits{M.D.}} \AND
\bauthor{\bsnm{Yukich}, \bfnm{J. E.}\binits{J.E.}}
(\byear{2005}).
\btitle{Normal approximation in geometric probability}.
In \bbooktitle{Stein's Method and Applications}.
\bseries{Lect. Notes Ser. Inst. Math. Sci. Natl. Univ. Singap.}
\bvolume{5}
\bpages{37--58}.
\blocation{Singapore}:
\bpublisher{Singapore Univ. Press}.
\bid{doi={10.1142/9789812567673_0003}, mr={2201885}}
\end{bincollection}
%
\iffalse\OrigBibText
Penrose, M.D. and Yukich, J.E. (2005) Normal approximation in geometric
probability. In Proceedings of the workshop `\textit{Stein's method
and applications}', Lecture notes series \textbf{5}, Singapore, pp. 37--58.
\endOrigBibText\fi
\bptok{imsref}%
% NOT OUTPUTTED:
% doi = http://dx.doi.org/10.1142/9789812567673_0003
\endbibitem

%b25 ###
\bibitem{RSZ}
%
\begin{barticle}[mr]
\bauthor{\bsnm{Reitzner}, \bfnm{M.}\binits{M.}},
\bauthor{\bsnm{Spodarev}, \bfnm{E.}\binits{E.}} \AND
\bauthor{\bsnm{Zaporozhets}, \bfnm{D.}\binits{D.}}
(\byear{2012}).
\btitle{Set reconstruction by {V}oronoi cells}.
\bjournal{Adv. in Appl. Probab.}
\bvolume{44}
\bpages{938--953}.
\bid{issn={0001-8678}, mr={3052844}}
\end{barticle}
%
\iffalse\OrigBibText
Reitzner, M., Spodarev, E. and Zaporozhets, D. (2012) Set
reconstruction by Voronoi cells. \textit{Adv. in Appl. Probab.}
\textbf{44}, 938--953.
\endOrigBibText\fi
\bptok{imsref}%
% NOT OUTPUTTED:
% url = http://projecteuclid.org/euclid.aap/1354716584
% number = 4
% fjournal = Advances in Applied Probability
\endbibitem

%b26 ###
\bibitem{SW}
%
\begin{bbook}[mr]
\bauthor{\bsnm{Schneider}, \bfnm{Rolf}\binits{R.}} \AND
\bauthor{\bsnm{Weil}, \bfnm{Wolfgang}\binits{W.}}
(\byear{2008}).
\btitle{Stochastic and Integral Geometry}.
\bseries{Probability and Its Applications}.
\blocation{Berlin}:
\bpublisher{Springer}.
\bid{doi={10.1007/978-3-540-78859-1}, mr={2455326}}
\end{bbook}
%
\iffalse\OrigBibText
Schneider, R. and Weil, W. (2008) \textit{Stochastic and Integral
Geometry}. Springer, Berlin.
\endOrigBibText\fi
\bptok{imsref}%
% NOT OUTPUTTED:
% doi = http://dx.doi.org/10.1007/978-3-540-78859-1
% isbn = 978-3-540-78858-4
% fpage = xii+693
\endbibitem

%b27 ###
\bibitem{Schreiber}
%
\begin{bincollection}[mr]
\bauthor{\bsnm{Schreiber}, \bfnm{Tomasz}\binits{T.}}
(\byear{2010}).
\btitle{Limit theorems in stochastic geometry}.
In \bbooktitle{New Perspectives in Stochastic Geometry}
\bpages{111--144}.
\blocation{Oxford}:
\bpublisher{Oxford Univ. Press}.
\bid{mr={2654677}}
\end{bincollection}
%
\iffalse\OrigBibText
Schreiber, T. (2010) Limit theorems in stochastic geometry. In `\textit
{New Perspectives in Stochastic Geometry}', pp. 111--144, Cambridge
University Press, Cambridge.
\endOrigBibText\fi
\bptok{imsref}%
\endbibitem

%b28 ###
\bibitem{SchulteVoronoi}
%
\begin{barticle}[mr]
\bauthor{\bsnm{Schulte}, \bfnm{Matthias}\binits{M.}}
(\byear{2012}).
\btitle{A central limit theorem for the {P}oisson--{V}oronoi approximation}.
\bjournal{Adv. in Appl. Math.}
\bvolume{49}
\bpages{285--306}.
\bid{doi={10.1016/j.aam.2012.08.001}, issn={0196-8858}, mr={3017961}}
\end{barticle}
%
\iffalse\OrigBibText
Schulte, M. (2012) A central limit theorem for the Poisson-Voronoi
approximation. \textit{Adv. in Appl. Math.} \textbf{49}, 285--306.
\endOrigBibText\fi
\bptok{imsref}%
% NOT OUTPUTTED:
% number = 3-5
% doi = http://dx.doi.org/10.1016/j.aam.2012.08.001
% fjournal = Advances in Applied Mathematics
\endbibitem

%b29 ###
\bibitem{TchZuyev}
%
\begin{barticle}[mr]
\bauthor{\bsnm{Tchoumatchenko}, \bfnm{Konstantin}\binits{K.}} \AND
\bauthor{\bsnm{Zuyev}, \bfnm{Sergei}\binits{S.}}
(\byear{2001}).
\btitle{Aggregate and fractal tessellations}.
\bjournal{Probab. Theory Related Fields}
\bvolume{121}
\bpages{198--218}.
\bid{doi={10.1007/PL00008802}, issn={0178-8051}, mr={1865485}}
\end{barticle}
%
\iffalse\OrigBibText
Tchoumatchenko, K. and Zuyev, S. (2001) Aggregate and fractal
tessellations. \textit{Probab. Th. Related Fields} \textbf{121}, 198--218.
\endOrigBibText\fi
\bptok{imsref}%
% NOT OUTPUTTED:
% number = 2
% doi = http://dx.doi.org/10.1007/PL00008802
% coden = PTRFEU
% fjournal = Probability Theory and Related Fields
\endbibitem

%b30 ###
\bibitem{YuLMN}
%
\begin{bincollection}[mr]
\bauthor{\bsnm{Yukich}, \bfnm{Joseph}\binits{J.}}
(\byear{2013}).
\btitle{Limit theorems in discrete stochastic geometry}.
In \bbooktitle{Stochastic Geometry, Spatial Statistics and Random Fields}.
\bseries{Lecture Notes in Math.}
\bvolume{2068}
\bpages{239--275}.
\blocation{Heidelberg}: \bpublisher{Springer}.
\bid{doi={10.1007/978-3-642-33305-7_8}, mr={3059650}}
\bptnote{check pages}%
\end{bincollection}
%
\iffalse\OrigBibText
Yukich, J.E. (2013) Limit theorems in discrete stochastic geometry. In
`\textit{Stochastic Geometry, Spatial Statistics and Random Fields}',
Springer, Berlin.
\endOrigBibText\fi
\bptok{imsref}%
% NOT OUTPUTTED:
% doi = http://dx.doi.org/10.1007/978-3-642-33305-7_8
\endbibitem

%b31 ###
\bibitem{Yu}
%
\begin{barticle}[mr]
\bauthor{\bsnm{Yukich}, \bfnm{J. E.}\binits{J.E.}}
(\byear{2015}).
\btitle{Surface order scaling in stochastic geometry}.
\bjournal{Ann. Appl. Probab.}
\bvolume{25}
\bpages{177--210}.
\bid{doi={10.1214/13-AAP992}, issn={1050-5164}, mr={3297770}}
\end{barticle}
%
\iffalse\OrigBibText
Yukich, J.E. (2015) Surface order scaling in stochastic geometry.
\textit{Ann. Appl. Probab.} \textbf{25}, 177--210.
\endOrigBibText\fi
\bptok{imsref}%
% NOT OUTPUTTED:
% number = 1
% doi = http://dx.doi.org/10.1214/13-AAP992
% fjournal = The Annals of Applied Probability
\endbibitem
\end{thebibliography}
\end{document}